\documentclass[12pt,a4paper]{article}

\usepackage[dvips]{color} 

\usepackage{amscd,amsmath,amssymb,amsfonts,times}
\usepackage{mathrsfs}
\usepackage[all]{xy}
\numberwithin{equation}{section}
    \newtheorem{thm}{Theorem}[section]
    \newtheorem{lem}[thm]{Lemma}
    
    \newtheorem{prop}[thm]{Proposition}
    \newtheorem{cor}[thm]{Corollary}

    \newtheorem{exmp}[thm]{Example}
    \newtheorem{rem}[thm]{Remark}
\newcommand{\qed}
{\mbox{}\nolinebreak$\square$\medbreak\par}
\newenvironment{pf}{\par\smallskip\noindent\emph{Proof.}}{\hfill\qed\par\smallskip}
\newenvironment{pf*}[1]{\par\smallskip\noindent\emph{#1.}}{\hfill\qed\par\smallskip}
\begin{document}
\title{A formula for Beilinson's regulator map
on $K_1$ of a fibration of curves having a totally degenerate 
semistable fiber}
\author{M. Asakura}
\date\empty
\maketitle
\tableofcontents

\def\can{\omega^*}
\def\canh{\omega}
\def\cano{\mathrm{canonical}}
\def\ff{{\mathit{false}}}
\def\Coker{\mathrm{Coker}}
\def\crys{\mathrm{crys}}
\def\zar{\mathrm{zar}}
\def\dlog{d{\mathrm{log}}}
\def\dR{{\mathrm{d\hspace{-0.2pt}R}}}            % de Rham
\def\et{{\mathrm{\acute{e}t}}}  % etale
\def\Frac{{\mathrm{Frac}}}
\def\phami{\phantom{-}}
\def\id{{\mathrm{id}}}              % identity
\def\Image{{\mathrm{Im}}}        % image
\def\Hom{{\mathrm{Hom}}}  
\def\ker{{\mathrm{Ker}}}          % kernel
\def\Pic{{\mathrm{Pic}}}
\def\CH{{\mathrm{CH}}}
\def\NS{{\mathrm{NS}}}
\def\NF{{\mathrm{NF}}}
\def\End{{\mathrm{End}}}
\def\pr{{\mathrm{pr}}}
\def\Proj{{\mathrm{Proj}}}
\def\ord{{\mathrm{ord}}}
\def\qis{{\mathrm{qis}}}
\def\reg{{\mathrm{reg}}}          %
\def\res{{\mathrm{res}}}          %
\def\Res{\mathrm{Res}}
\def\Spec{{\mathrm{Spec}}}     % spectrum
\def\syn{{\mathrm{syn}}}
\def\cont{{\mathrm{cont}}}
\def\ind{{\mathrm{ind}}}
\def\inv{{\mathrm{inv}}}
\def\dec{{\mathrm{dec}}}
\def\Ext{{\mathrm{Ext}}}
\def\MHS{{\mathrm{MHS}}}
\def\Gr{{\mathrm{Gr}}}
\def\cand{{\mathrm{ex}}}

\def\bA{{\mathbb A}}
\def\bC{{\mathbb C}}
\def\C{{\mathbb C}}
\def\G{{\mathbb G}}
\def\bE{{\mathbb E}}
\def\bF{{\mathbb F}}
\def\F{{\mathbb F}}
\def\bH{{\mathbb H}}
\def\bJ{{\mathbb J}}
\def\bN{{\mathbb N}}
\def\bP{{\mathbb P}}
\def\P{{\mathbb P}}
\def\bQ{{\mathbb Q}}
\def\Q{{\mathbb Q}}
\def\bR{{\mathbb R}}
\def\R{{\mathbb R}}
\def\bZ{{\mathbb Z}}
\def\Z{{\mathbb Z}}
\def\cA{{\mathscr A}}
\def\cD{{\mathscr D}}
\def\cM{{\mathscr M}}
\def\cL{{\mathscr L}}
\def\cE{{\mathscr E}}
\def\cO{{\mathscr O}}
\def\O{{\mathscr O}}
\def\cR{{\mathscr R}}
\def\cS{{\mathscr S}}
\def\cX{{\mathscr X}}
\def\cH{{\mathscr H}}
\def\ccH{{\mathscr H}_C}
\def\PF{{\mathit{PF}}}
\def\Div{{\mathrm{Div}}}
\def\codim{{\mathrm{codim}}}
\def\Van{{\mathrm{Ev}_P}}
%                                 Greece
%
\def\ep{\epsilon}
\def\vG{\varGamma}
\def\vg{\varGamma}

%
%                                 simple
%
%
%
\def\lra{\longrightarrow}
\def\lla{\longleftarrow}
\def\Lra{\Longrightarrow}
\def\hra{\hookrightarrow}
\def\lmt{\longmapsto}
\def\ot{\otimes}
\def\op{\oplus}
%                              decolation
\def\wt#1{\widetilde{#1}}
\def\wh#1{\widehat{#1}}
\def\spt{\sptilde}
\def\ol#1{\overline{#1}}
\def\ul#1{\underline{#1}}
\def\us#1#2{\underset{#1}{#2}}
\def\os#1#2{\overset{#1}{#2}}
\def\lim#1{\us{#1}{\varinjlim}}
\def\plim#1{\us{#1}{\varprojlim}}

\section{Introduction}
Let $X$ be a projective nonsingular
variety over the complex number field $\C$.
Let $H^i_\cM(X,\Z(j))$ denotes the motivic cohomology group. 
It is known that 
$H^i_\cM(X,\Q(j))$ is isomorphic to Quillen's $K$-group $K_{2j-i}(X)^{(j)}$.
By the theory of higher Chern classes,
we have the {\it Beilinson regulator map} (higher Chern class map)
\[
\reg_{i,j}:H^i_\cM(X,\Z(j))\lra H_\cD^i(X,\Z(j))
\]
to the Deligne-Beilinson cohomology group
(\cite{schneider}).
The purpose of this paper is to provide a systematic method for computations
of the regulator map in case that $(i,j)=(3,2)$ (namely $K_1$)
and $X$ a fibration of curves having a totally degenerate semistable fiber.
Here we mean by a {\it fibration of curves} a surjective morphism $f:X\to C$ with a section 
$e:C\to X$ where $X$ (resp. $C$) is a projective smooth surface (resp. curve).
A fiber is called {\it semistable} if it is reduced and a
formal neighborhood of each singular point
is isomorphic to ``$xy=0$".
A fiber is called {\it totally degenerate} if each component is a rational curve.

\medskip

The cup-product pairing gives rise to a map 
$\C^\times\ot \Pic(X)\cong \C^\times\ot H^2_\cM(X,\Z(1))
\to H^3_\cM(X,\Z(2))$. Its image is called the {\it decomposable} part,
and the cokernel is called the {\it indecomposable} part.
The decomposable part does not affect serious difficulty,
while the indecomposable part plays the central role in the study of
$H^3_\cM(X,\Z(2))$.
According to \cite{lewisJAG}, we call an element $\xi\in H^3_\cM(X,\Z(2))$
{\it regulator indecomposable} if
$\reg_{3,2}(\xi)$ does not lie in the image of $\C^\times\ot\NS(X)$.
Obviously regulator indecomposable elements are indecomposable.
The converse is also true if the Beilinson-Hodge conjecture for $K_2$ is true.
Lewis and Gordon constructed regulator indecomposable elements in case
$X$ is a product of `general' elliptic curves (\cite{lewisJAG} Theorem 1).
There are lots of other related works, though I don't catch up all of them.
On the other hand, in case that $X$ is defined over a number field,
the question is more difficult, and as far as I know
there are only a few of such examples (e.g. \cite{R} \S 12).

In this paper we give a new method for computation of the composition of the maps
\[
H_{\cM,D}^3(X,\Q(2))\lra H_{\cM}^3(X,\Q(2))\os{\reg_{3,2}}{\lra}
H_\cD^3(X,\Q(2))
\]
where $D=f^{-1}(P)$ is a singular fiber, $H_{\cM,D}^3(X,\Q(2))$ denotes the motivic cohomology
supported on $D$ (Theorem \ref{ExtMHS}).
As an application we give a number of
examples of regulator indecomposable elements (\S \ref{Example-sect}).
One of the technical key points is to use
certain rational 2-forms ``$\Lambda(X)_\mathrm{rat}$" which will be introduced
in \S \ref{good-sect}. 
Our method has an advantage in explicit computations
and it also works in case that the base field
is a number field.
I hope that it will bring a new progress in the study of
 Beilinson's conjecture on special values of $L$-functions (\cite{RSC}), though
I don't have a result in this direction so far.

\medskip

This paper is organized as follows.
\S \ref{reg-sect} is a quick review of $H^3_\cM(X,\Q(2))$ and Beilinson regulator.
\S \ref{deRh-sect} is a preparation to state and prove the main theorem 
(Theorem \ref{ExtMHS}) in \S \ref{exp-sect}.
In \S \ref{Example-sect} we give examples 
of regulator indecomposable elements for certain elliptic surfaces
defined over $\Q$ with arbitrary large $h^{2,0}$ (Cor.\ref{Example-cor}).
\S \ref{Appendix-sect} is an appendix providing how to compute 
the Gauss-Manin connection for a hyper elliptic fibrations.

\bigskip

\noindent{\bf Acknowledgment.}
A rough idea was inspired
during my visit to the University of Alberta in September 2012,
especially when I discussed the paper \cite{lewis}
 with Professor James Lewis.
I would like to express special thanks to him.
I'd also like to thank the university members
for their hospitality.

%%%%%%%%%%%%%%%%%%%%%%%%%%%%%%%%%%%%%%%%%%%%%%%
\section{Real Regulator map on $H^3_\cM(X,\Q(2))$}\label{reg-sect}
For a variety $X$ over a field $K$ of characteristic zero,
we denote by $H^\bullet_\dR(X)=H^\bullet_\dR(X/K)$ (resp.
$H_\bullet^\dR(X)=H_\bullet^\dR(X/K)$) 
the de Rham cohomology (resp. de Rham homology) cf. \cite{h}.
When $K=\C$
we denote by $H^\bullet_B(X,\Q)=H^\bullet_B(X(\C),\Q)$ 
(resp. $H_\bullet^B(X,\Q)$) the Betti cohomology (resp. Betti homology).

\subsection{$H^3_\cM(X,\Q(2))$ and indecomposable parts}
Let $H_{\cM,Z}^i(X,\Q(j))$ denotes the motivic cohomology of
a smooth variety $X$ supported on a closed subscheme $Z\subset X$.
They fit into the localization exact sequence
\[
\cdots\lra
H^i_{\cM,Z}(X,\Q(j))\lra H^i_{\cM}(X,\Q(j))\lra H^i_{\cM}(X\setminus Z,\Q(j))
\lra \cdots.
\]
Of particular interest to us is the case $(i,j)=(3,2)$.
Let us describe $H^3_\cM(X,\Q(2))$ explicitly.
Let $Z_i(X)=Z^{\dim X-i}(X)$ be the free abelian group of irreducible subvarieties
of dimension $i$.
We denote by $\eta_V$ the field of rational functions on an integral scheme $V$.
Let $Z\subset X$ be an irreducible divisor, and $\wt{Z}\to Z$
the normalization.
Let $j:\wt{Z}\to Z\hra X$ be the composition.
Then 
we define $\Div_Z(f):=j_*\Div_{\wt{Z}}(f)\in Z^2(X)$ the push-forward of the Weil divisor
on $\wt{Z}$ by $j$.
Let
\[
\partial_1:\bigoplus_{\codim Z=1} \eta_Z^\times\lra Z^2(X),\quad
[f,Z]\longmapsto \Div_Z(f)
\]
be a homomorphism where we write
\[
[f,Z]:=(\cdots,1,f,1,\cdots)\in \bigoplus_{\codim Z=1} \eta_Z^\times,\quad
\mbox{($f$ is placed in the $Z$-component).}
\] 
Let $Z=\sum_{i=1}^r Z_i$ be a dvisor with $Z_i$ irreducible.
Then, for a smooth open set $Z^\circ\subset Z$ we have
\begin{equation}\label{mot-0-1}
H^3_{\cM,Z}(X,\Q(2))\os{\text{can.}}{\lra}
H^3_{\cM,Z^\circ}(X^\circ,\Q(2))\os{\cong}{\longleftarrow} \O(Z^\circ)^\times\ot\Q
\end{equation}
where $X^\circ:=X\setminus(Z\setminus Z^\circ)$, and this induces a canonical isomorphism
\begin{equation}\label{mot-0-0}
H^3_{\cM,Z}(X,\Q(2))\cong \ker[\bigoplus_{i=1}^r \eta^\times_{Z_i}
\os{\partial_1}{\lra} Z^2(X)]\ot\Q.
\end{equation}
Let \[\partial_2:K^M_2(\eta_X)\lra \bigoplus_{\codim Z=1} \eta_Z^\times\]
\[
\partial_2\{f,g\}=\sum_{\codim Z=1}
\left[(-1)^{\ord_Z(f)\ord_Z(g)}\frac{f^{\ord_Z(g)}}{g^{\ord_Z(f)}}|_Z,Z\right]
\]
be the tame symbol.
Then the natural map 
$H^3_{\cM,Z}(X,\Q(2))\to H^3_{\cM}(X,\Q(2))$ together with \eqref{mot-0-0}
induces the isomorphism
\begin{equation}\label{mot-0}
H^3_\cM(X,\Q(2))\cong \left(\frac{\ker(\bigoplus \eta_Z^\times
\os{\partial_1}{\lra} Z^2(X))}
{\Image(K^M_2(\eta_X)\os{\partial_2}{\lra} \bigoplus \eta_Z^\times)}\right)\ot\Q.
\end{equation}

Let $K$ be the base field of $X$ and $L/K$ a finite extension. Write $X_L:=X\times_KL$.
Then there is the obvious map
\[
L^\times\ot Z^1(X_L)\lra H^3_\cM(X_L,\Q(2)),\quad \lambda\ot D\longmapsto[\lambda,D].
\]
Let $N_{L/K}:H^3_\cM(X_L,\Q(2))\to H^3_\cM(X,\Q(2))$ be the norm map on motivic 
cohomology.
Then we put
\[
H^3_\cM(X,\Q(2))_\dec:=\sum_{[L:K]<\infty}
N_{L/K}(\Image(L^\times\ot Z^1(X_L)\to H^3_\cM(X_L,\Q(2))))
\]
and call it the {\it decomposable part}.
We put
\[
H^3_\cM(X,\Q(2))_\ind:=
H^3_\cM(X,\Q(2))/H^3_\cM(X,\Q(2))_\dec
\]
and call it the {\it indecomposable part}.
The indecomposable part plays the central role in the study of $H^3_\cM(X,\Q(2))$.

\subsection{Beilinson regulator on indecomposable parts}\label{bregin-sect}
Let $X$ be a smooth projective variety over $\C$.
By the theory of universal Chern class, there is the {\it Beilinson regulator map}
\begin{equation}\label{beireg-0}
\reg:H^3_{\cM}(X,\Q(2))\lra H^3_{\cD}(X,\Q(2))
\cong\Ext_\MHS^1(\Q,H^2(X,\Q(2)))
\end{equation}
to the Deligne-Beilinson cohomology group, which is isomorphic to
the Yoneda extension group of mixed Hodge structures
where $H^2(X,\Q(2))=(H^2_B(X,\Q(2)),F^\bullet H^2_\dR(X))$ denotes the Hodge structure
(of weight $-2$).
We here write down
the regulator map \eqref{beireg-0} in terms of extension of mixed Hodge structure.

Let $Z=\sum_{i=1}^r Z_i\subset X$ be a divisor.
There is also the Beilinson regulator map on $H^3_{\cM,Z}(X,\Q(2))$, which we denote by
$\reg_Z$.
Let us consider a commutative diagram
\begin{equation}\label{beireg-0-3}
\xymatrix{
H^3_{\cM}(X,\Q(2))\ar[d]_{\reg}& H^3_{\cM,Z}(X,\Q(2))\ar[l]\ar@{^{(}->}[r]
\ar[d]_{\reg_Z}& H^3_{\cM,Z^\circ}(X^\circ,\Q(2))\ar[d]_{\reg_{Z^\circ}}
&\ar[l]_{\quad\cong}^{\quad\eqref{mot-0-1}} \ar[d]^{\mathrm{dlog}}\O(Z^\circ)^\times\ot\Q\\
H^3_{\cD}(X,\Q(2))& H^3_{\cD,Z}(X,\Q(2))\ar[l]_\iota\ar@{^{(}->}[r]
& H^3_{\cD,Z^\circ}(X^\circ,\Q(2))
&\ar[l]_{\quad\cong} H^1_\cD(Z^\circ,\Q(1)).
}
\end{equation}
Here the commutativity of the right square follows from the
Riemann-Roch theorem without denominator (\cite{gillet}), and the others follow from
the functoriality of the regulator map.
The middle and right squares in \eqref{beireg-0-3} and the isomorphism
\eqref{mot-0-1} induce a map
\begin{equation}\label{beireg-0-1}
\ker[\bigoplus_{i=1}^r\eta_{Z_i}^\times\os{\partial_1}{\to} Z^2(X)]
\lra H^3_Z(X,\Q(2))\cap H^{0,0},\quad
\xi\longmapsto\nu_\xi,
\end{equation}
where we write $M\cap H^{0,0}:=\Hom_\MHS(\Q,M)$.
This is characterized by 
\[
\nu_\xi|_{Z^\circ}=\left(\frac{df_i}{f_i}\right)\in 
H^3_{Z^\circ}(X^\circ,\Q(2))\cong H^1(Z^\circ,\Z(1)),\quad \xi=\sum_{i=1}^r [f_i,Z_i].
\]
Let $\langle Z_i\rangle\subset H^2(X,\Q(2))$ denotes the subgroup
generated by the cycle classes of $Z_i$.
Then the map $\iota$ in \eqref{beireg-0-3} induces a commutative diagram
\begin{equation}\label{beireg-0-2}
\xymatrix{
H^3_{\cD,Z}(X,\Q(2))\ar[r]\ar[d]_\iota& H^3_{Z}(X,\Q(2))\cap H^{0,0}\ar[d]^\delta\\
H^3_{\cD}(X,\Q(2))\ar[r]& \Ext_\MHS(\Q,H^2(X,\Q(2))/\langle Z_i\rangle)\\
}\end{equation}
where $\delta$ is the connecting homomorphism
arising from the localization exact sequence
\[
\cdots\lra H^j_Z(X,\Q(2))\lra H^j(X,\Q(2))
\lra
H^j(X-Z,\Q(2))\lra \cdots.
\]
Summing up the above, we have the following theorem. 
\begin{thm}[e.g. \cite{sato} 11.2]\label{bregin-thm}
The following diagram is commutative
\[
\xymatrix{
\ker[\bigoplus_{i=1}^r\eta_{Z_i}^\times\os{\partial_1}{\to} Z^2(X)]\ar[rd]^{\eqref{beireg-0-1}}
\ar[d]_\cong^{\eqref{mot-0-0}}\\
H^3_{\cM,Z}(X,\Q(2))\ar[d]\ar[r]^{\reg_Z\qquad}&H_Z^3(X,\Q(2))\cap H^{0,0}
\ar[d]^{\delta~\eqref{beireg-0-2}}\\
H^3_{\cM}(X,\Q(2))\ar[r]^{\reg\hspace{1.5cm}}&
\Ext_\MHS^1(\Q,H^2(X,\Q(2))/\langle Z_i\rangle).
}
\]
\end{thm}
For the later use, we write down $\delta$ more explicitly.
Let $n=\dim X$.
Under the isomorphism of Poincare-Lefschetz duality, the map 
$\delta$ coincides with the connecting homomorphism arising from the exact sequence
\[
\cdots\lra H_i(X,\Q(2-n))\lra H_i(X,Z;\Q(2-n))\lra
H_i(Z,\Q(2-n))\lra \cdots.
\]
Put $H^i_\dR(X)_Z:=\ker[ H^i_\dR(X)\to H^i_\dR(Z)]$
and 
\[
M:=H_{2n-2}(X,\Q(2-n))/H_{2n-2}(Z,\Q(2-n))\cong 
H^2(X,\Q(2))/\langle Z_i\rangle.
\]
Then there is the natural isomorphism
$M\ot_\Q\C\cong
\Hom(H^{2n-2}_\dR(X)_Z,\C)
$
and it induces
\begin{equation}
\Ext^1_{\MHS}(\Q,M)\cong
\Coker[
H_{2n-2}(X,\Q(2-n))\os{\Phi}{\to} \Hom(F^{n-1}H^{2n-2}_\dR(X)_Z,\C)
]\label{beireg-2}
\end{equation}
where $\Phi$ is defined as follows
\[
\Phi(\Delta)=\left[
\omega\longmapsto \int_\Delta\omega
\right],\quad \omega\in F^{n-1}H^{2n-2}_\dR(X).
\]
Let $\omega_{X,Z}\in F^{n-1}H^{2n-2}_\dR(X,Z)$ denotes a lifting of 
$\omega\in F^{n-1}H^{2n-2}_\dR(X)_Z$ via the surjective map
$H^{2n-2}_\dR(X,Z)\to H^{2n-2}_\dR(X)_Z$.
Let $\gamma\in H_{2n-3}(Z,\Q(2-n))\cap H^{0,0}
\cong H^3_Z(X,\Q(2))\cap H^{0,0}$.
Let $\Gamma\in H_{2n-2}(X,Z;\Q(2-n))$ be an arbitrary element 
such that $\partial(\Gamma)=\gamma$.
Then we have
\begin{equation}\label{beireg-3}
\delta(\gamma)=\left[
\omega\longmapsto \int_\Gamma\omega_{X,D}
\right]
\end{equation}
under the isomorphism \eqref{beireg-2}.

\begin{prop}\label{bei-hodge} 
Let $Z=\sum_{i=1}^r Z_i$ be a divisor. Then
\[
\reg_Z:H^3_{\cM,Z}(X,\Q(2))\lra H^3_Z(X,\Q(2))\cap H^{0,0}
\cong H_{2n-3}(Z,\Q(2-n))\cap H^{0,0}
\]
is surjective.
\end{prop}
\begin{pf}
Let $Z^\circ=\cup Z_i^\circ\subset Z$ be a regular locus, and put
$\Sigma:=Z\setminus Z^\circ$ and $X^\circ:=X\setminus \Sigma$.
Let $\Q\Sigma\subset Z^2(X)$ be the subgroup generated by components of $\Sigma$
of codimension 2. 
Then we have a commutative diagram
\[
\xymatrix{
0\ar[r]&H^3_{\cM,Z}(X,\Q(2))\ar[d]_{\reg_Z}\ar[r]&
\O(Z^\circ)^\times\ot\Q\ar[r]^{\quad \partial_1}\ar[d]^{\mathrm{dlog}}
&\Q\Sigma\ar@{^{(}->}[d]\\
0\ar[r]&H^3_Z(X,\Q(2))\cap H^{0,0}\ar[r]&
H^3_{Z^\circ}(X^\circ,\Q(2))\cap H^{0,0}\ar[r]
&H^4_\Sigma(X,\Q(2))\\
&&H^1(Z^\circ,\Q(1))\cap H^{0,0}\ar[u]_\cong
}
\]
with exact rows. Now the assertion follows from the surjectivity of $\mathrm{dlog}$.
\end{pf}

Put
\[
H^2_B(X)_\ind:=H^2_B(X,\Q(1))/\NS(X)\ot\Q,\quad
H^2_\dR(X)_\ind:=H^2_\dR(X)/\NS(X)\ot\C,
\]
\[
H^2(X)_\ind:=(H^2_B(X)_\ind,F^\bullet H^2_\dR(X)_\ind)
\mbox{ (= a Hodge structure of weight 0).}
\]
Then the Beilinson regulator map \eqref{beireg-0} yields a commutative diagram
\begin{equation}\label{bregin-2}
\xymatrix{
0\ar[d]&0\ar[d]\\
H^3_\cM(X,\Q(2))_\dec\ar[d]\ar[r]&\Ext_\MHS^1(\Q,\NS(X)\ot\Q(1)))\ar[d]\\
H^3_\cM(X,\Q(2))\ar[d]\ar[r]^{\reg\qquad}&\Ext^1_\MHS(\Q,H^2(X,\Q(2)))\ar[d]\\
H^3_\cM(X,\Q(2))_\ind\ar[d]\ar[r]
&\Ext_\MHS^1(\Q,H^2(X)_\ind\ot\Q(1))\ar[d]\\
0&0
}
\end{equation}
The top arrow is simply written by ``$\log$", namely the composition
\[
\C^\times\ot\Pic(X)\to H^3_\cM(X,\Q(2))_\dec\lra 
\Ext_\MHS^1(\Q,\NS(X)\ot\Q(1)))\cong \C/\Q(1)\ot\NS(X)
\]
is given by $\lambda\ot Z\mapsto \log(\lambda)\ot Z$.
The bottom arrow, the {\it regulator map on indecomposable parts},
is the main research subject of this paper.

The {\it real regulator map} is the composition of $\reg$ and the canonical map
\[
\Ext^1_\MHS(\Q,H^2(X,\Q(2)))\lra
\Ext^1_{\R\mbox{-}\MHS}(\R,H^2(X,\R(2)))
\]
to the extension group of real mixed Hodge structures, 
which we denote by $\reg_\R$:
\begin{equation}\label{bregin-3}
\reg_\R:H^3_\cM(X,\Q(2))\lra \Ext_{\R\mbox{-}\MHS}^1(\R,H^2(X,\R(2)))
\cong H^2_B(X,\R(1))\cap H^{1,1}.
\end{equation}
This also induces the map
\begin{equation}
\reg_\R:H^3_\cM(X,\Q(2))_\ind\to \Ext_{\R\mbox{-}\MHS}^1(\R,H^2(X)_\ind\ot\R(1))
\cong (H^2_B(X)_\ind\ot\R)\cap H^{1,1}
\end{equation}
on the indecomposable part (we use the same symbol since there will not be confusion).

\subsection{$\Q$-structure on determinant of $H^3_\cD(X/\R,\R(2))$}
Suppose that $X$ is a projective smooth variety over $\Q$.
Write $X_\C:=X\times_\Q\C$.
The {\it infinite Frobenius} map $F_\infty$ is defined to be the
anti-holomorphic map on $X(\C)=\mathrm{Mor}_\Q(\Spec\C,X)$ 
induced from the complex conjugation on $\Spec\C$.
For a subring $A\subset \R$,
the infinite Frobenius map 
acts on the Deligne-Beilinson complex $A_X(j)_\cD$ in a canonical way,
so that we have the involution on $H^\bullet_\cD(X_\C,A(j))$, which we denote by the
same notation $F_\infty$.
We define
\[
H^\bullet_\cD(X/\R,A(j)):=H^\bullet_\cD(X_\C,A(j))^{F_\infty=1}
\]
the fixed part by $F_\infty$. We call it the {\it real Deligne-Beilinson cohomology}.
Since the action of $F_\infty$ is compatible via the Beilinson regulator map,
we have
\begin{align}
\reg_\R:H^3_\cM(X,\Q(2))\lra H_{\cD}:=&
\Ext_{\R\mbox{-}\MHS}^1(\R,H^2(X_\C,\R(2)))^{F_\infty=1}\\
\cong& \frac{H^2_B(X_\C,\R(1))^{F_\infty=1}}{F^2H^2_\dR(X/\R)},\label{qstr-3}
\end{align}
and 
\begin{align}
\reg_\R:H^3_\cM(X,\Q(2))\lra H_{\cD,\ind}:=&
\Ext_{\R\mbox{-}\MHS}^1(\R,H^2(X_\C)_\ind\ot\R(1))^{F_\infty=1}\\
\cong& \frac{[H^2_B(X_\C)_\ind\ot\R]^{F_\infty=1}}{F^2H^2_\dR(X/\R)}.\label{qstr-4}
\end{align}
%(see \eqref{bregin-1} for the notation ``$H^2_B(X)_\ind$"). 
There are the canonical $\Q$-structures $e_\Q$ and $e_{\ind,\Q}$ 
on the determinant vector spaces $\det H_\cD$ and
$\det H_{\cD,\ind}$:
\[
\R\cdot e_\Q=\det H_\cD,\quad
\R\cdot e_{\ind,\Q}=\det H_{\cD,\ind}.
\]
Here we recall the definition.
The isomorphisms \eqref{qstr-3} and \eqref{qstr-4} induce 
\begin{equation}\label{qstr-1}
\det H_\cD
\cong 
\det [H^2_B(X_\C,\R(1))^{F_\infty=1}]\ot[\det F^2H^2_\dR(X/\R)]^{-1},
\end{equation}
and 
\begin{equation}\label{qstr-2}
\det H_{\cD,\ind}
\cong \det[(H^2_B(X_\C)_\ind\ot\R)^{F_\infty=1}]
\ot[\det F^2H^2_\dR(X/\R)]^{-1}.
\end{equation}
The right hand sides of \eqref{qstr-1} and \eqref{qstr-2} have the $\Q$-structures
induced from the $\Q$-structures
\[
H^2_B(X_\C,\Q(1))^{F_\infty=1},\quad
H^2_B(X_\C)_\ind^{F_\infty=1}, \quad F^2H^2_\dR(X/\Q).
\]
The $\Q$-structures $e_\Q$ and $e_{\ind,\Q}$ are defined to be the corresponding
one:
\begin{equation}\label{qstr-9}
\Q\cdot e_\Q\cong 
\det [H^2_B(X_\C,\Q(1))^{F_\infty=1}]\ot[\det F^2H^2_\dR(X/\Q)]^{-1},
\end{equation}
\begin{equation}\label{qstr-10}
\Q\cdot e_{\Q,\ind}\cong 
\det[H^2_B(X_\C)_\ind^{F_\infty=1}]
\ot[\det F^2H^2_\dR(X/\Q)]^{-1}.
\end{equation}
\subsection{$e_\Q^\ff$ and $e_{\ind,\Q}^\ff$}\label{false-sect}
We introduce other $\Q$-structures $e_\Q^\ff$ and $e_{\ind,\Q}^\ff$ on $\det H_\cD$ and
$\det H_{\cD,\ind}$. 
For simplicity, we assume $\dim X=2$.
Put
\[
H_2(X_\C,\Q)_\ind:=H_2(X_\C,\Q)/(\NS(X_\C)\ot\Q(1))\cong 
H^2_B(X_\C)_\ind\ot\Q(1),
\]
\[
H^2_\dR(X/\Q)_\ind:=\mathrm{Coim}(H^2_\dR(X/\Q)\lra H^2_\dR(X/\C)/(\NS(X_\C)\ot\C)).
\]
Note that $H^2_\dR(X/\Q)_\ind\ot\C\os{\cong}{\to}
H^2_\dR(X/\C)/(\NS(X_\C)\ot\C)$.
There are exact sequences
\begin{equation}\label{qstr-5}
0\lra H_2(X_\C,\R)^{F_\infty=1}\lra \Hom(F^1H^2_\dR(X/\Q),\R)\lra H_\cD\lra 0
\end{equation}
\begin{equation}\label{qstr-6}
0\lra H_2(X_\C,\R)_\ind^{F_\infty=1}\lra \Hom(F^1H^2_\dR(X/\Q)_\ind,\R)\lra 
H_{\cD,\ind}\lra 0
\end{equation}
under the canonical isomorphisms
\begin{equation}\label{qstr-7}
H^2_B(X_\C,\C)\cong H^2_\dR(X/\C),\quad H^2_B(X_\C,\Q(2))\cong H_2(X_\C,\Q).
\end{equation}
Then we define $e^{\ff}_\Q$ and $e^\ff_{\ind,\Q}$ the $\Q$-structures induced from
\[
H_2(X_\C,\Q)^{F_\infty=1},\quad H_2(X_\C,\R)_\ind^{F_\infty=1},\quad
H^2_\dR(X/\Q),\quad
H^2_\dR(X/\Q)_\ind.
\]
Hence we have
\begin{equation}\label{qstr-11}
\Q\cdot e_\Q^\ff\cong 
[\det H_2(X_\C,\Q)^{F_\infty=1}]^{-1}\ot[\det F^1H^2_\dR(X/\Q)]^{-1},
\end{equation}
\begin{equation}\label{qstr-12}
\Q\cdot e_{\ind,\Q}^\ff\cong 
[\det H_2(X_\C,\Q)_\ind^{F_\infty=1}]^{-1}\ot[\det F^1H^2_\dR(X/\Q)_\ind]^{-1}.
\end{equation}

\begin{prop}\label{qstr-8}
Put\[
r:=\dim H_2(X_\C,\Q)^{F_\infty=1}=\dim H^2_B(X_\C,\Q(1))^{F_\infty=-1},
\] 
\[
s:=\dim H_2(X_\C,\Q)_\ind^{F_\infty=1}=\dim H^2_B(X_\C)_\ind^{F_\infty=-1}
=r-\dim \NS(X_\C)^{F_\infty=-1}.
\]
Write \[
H_B:=H^2_B(X_\C,\Q(1)), \quad H_{B,\ind}:=H^2_B(X_\C)_\ind,\] 
\[
F^\bullet H_\dR:=F^\bullet H^2_\dR(X/\Q),\quad
F^\bullet H_{\dR,\ind}:=F^\bullet H^2_\dR(X/\Q)_\ind\] simply.
Then
\[
\Q\cdot e^{\ff}_\Q=\Q\cdot e_\Q\ot\Q(-r)\ot \det H_\dR\ot
[\det H_B]^{-1},
\]
\[
\Q\cdot e^{\ff}_{\ind,\Q}=\Q\cdot e_{\ind,\Q}\ot\Q(-s)\ot \det H_{\dR,\ind}\ot
[\det H_{B,\ind}]^{-1},
\]
where we mean
\[
\det H_\dR\ot
[\det H_B]^{-1}\subset \det H_\dR^2(X/\C)\ot
[\det H_B^2(X_\C,\C)]^{-1}\os{\eqref{qstr-7}}{\cong}\C,\mbox{ etc.}\]
\end{prop}
\begin{pf}
By the Poincare duality, one has
\[
\det F^1H_\dR=[\det H_\dR]^{-1}\ot \det F^2H_\dR,\quad
\det F^1H_{\dR,\ind}=[\det H_{\dR,\ind}]^{-1}\ot \det F^2H_\dR.
\]
Moreover one has
\[
\det [H_2(X_\C,\Q)^{F_\infty=1}]
=\det [H^2_B(X_\C,\Q(2))^{F_\infty=1}]
=\det [H_B^{F_\infty=-1}\ot\Q(1)]
=\Q(r)\ot \det H_B^{F_\infty=-1}
\]
and
\[
\det [H_2(X_\C,\Q)_\ind]^{F_\infty=1}
=\det[ H_{B,\ind}^{F_\infty=-1}\ot\Q(1)]
=\Q(s)\ot\det H_{B,\ind}^{F_\infty=-1}.
\]
Therefore we have
\begin{align*}
\Q\cdot e^{\ff}_\Q\ot e_\Q^{-1}
&=
\Q(-r)\ot[\det H_B^{F_\infty=-1}]^{-1}
\ot [\det H_B^{F_\infty=1}]^{-1}\ot[\det H_\dR]\\
&=
\Q(-r)\ot[\det H_B]^{-1}
\ot[\det H_\dR]\\
\end{align*}
by \eqref{qstr-9} and \eqref{qstr-11},
and
\begin{align*}
\Q\cdot e^{\ff}_{\ind,\Q}\ot e_{\ind,\Q}^{-1}
&=
\Q(-s)\ot[\det H_{B,\ind}^{F_\infty=-1}]^{-1}
\ot [\det H_{B,\ind}^{F_\infty=1}]^{-1}\ot[\det H_{\dR,\ind}]\\
&=
\Q(-s)\ot[\det H_{B,\ind}]^{-1}
\ot[\det H_{\dR,\ind}]\\
\end{align*}
by \eqref{qstr-10} and \eqref{qstr-12}. This completes the proof.
\end{pf}
\begin{rem}
The Poincare duality implies
\[
(\det H_B)^{\ot 2}\cong 
 H^4_B(X_\C,\Q(2))^{\ot m}\cong
H^4_\dR(X/\Q)^{\ot m}\cong
(\det H_\dR)^{\ot 2},
\]
and
\[
(\det H_{B,\ind})^{\ot 2}\cong 
 H^4_B(X_\C,\Q(2))^{\ot m'}\cong
H^4_\dR(X/\Q)^{\ot m'}\cong
(\det H_{\dR,\ind})^{\ot 2}.
\]
Therefore $(\det H_\dR\ot[\det H_B]^{-1})$ and
$(\det H_{\dR,\ind}\ot[\det H_{B,\ind}]^{-1})$ are contained in $\sqrt{\Q^\times}$
(possibly rational numbers).
\end{rem}

%%%%%%%%%%%%%%%%%%%%%%%%%%%%%%%%%%%%%%%%%%%%%%%

\section{Cohomology of Fibration of curves and rational 2-forms}\label{deRh-sect}

\subsection{Notation}\label{deRh-sect-1}
Let $X$ (resp. $C$) is a projective smooth surface (resp. curve) over 
$K$ a field of characteristic 0, and
let $f:X\to C$ be a surjective morphism with a section $e:C\to X$.
The general fiber $X_t:=f^{-1}(t)$ is a projective smooth curve of genus $g>0$.
Throughout this section, we use the following notation.
\begin{itemize}
\item
Write $X_{\ol{K}}=X\times_K\ol{K}$. 
Define $\NF(X_{\ol{K}})\subset \NS(X_{\ol{K}})$ to be the subgroup of the Neron-Severi group
generated by the section $e(C)$ and fibral divisors (i.e. 
irreducible components of singular fibers).
\item
$\NF_\dR(X):=H^2_\dR(X)\cap (\NF(X_{\ol{K}})\ot_\Z\ol{K})\subset H^2_\dR(X_{\ol{K}}/\ol{K})$.
\item
For a Zariski open set $S\subset C$ and $V:=f^{-1}(S)$, we put
\[
H^2_\dR(V)_0:=\ker[
H^2_\dR(V)\lra\prod_{s\in S}H^2_\dR(f^{-1}(s))\times H^2_\dR(e(S))]
\]
where the arrow is the restriction map.
Note $H^2_\dR(e(S))=0$ unless $S=C$.
Note also that ``$f^{-1}(s)$" suffices to run over only singular fibers
and one smooth fiber.
\item
$\NF_\dR(V):=\Image[\NF_\dR(X)\to H^2_\dR(V)]$.
\end{itemize}
\begin{rem}\label{nf-rem-dege}
The intersection pairing $\NF(X_{\ol{K}})\ot \NF(X_{\ol{K}})\to \Q$ is non-degenerate.
This follows from Zariski's lemma (\cite{barth} III (8.2)).
\end{rem}
\begin{rem}\label{nf-rem}
$\NF_\dR(X)\ot_K\ol{K}=\NF(X_{\ol{K}})\ot_\Z\ol{K}
=\NF_\dR(X_{\ol{K}})$ in $H^2_\dR(X_{\ol{K}}/\ol{K})$, and
hence $\NF_\dR(V)\ot_K\ol{K}
=\NF_\dR(V_{\ol{K}})$.
This is proven by using \cite{AEC} II Lemma 5.8.1.
\end{rem}
\begin{rem}\label{vdr-nf-rem}
Let $\NF_\dR(X)^\perp$ denotes the
orthogonal complements of 
$\NF_\dR(X)$ in $H^2_\dR(X)$
with respect to the cup-product pairing.
Then
$\NF_\dR(X)^\perp= H^2_\dR(X)_0$ by definition, and hence
$\NF_\dR(X)\op H^2_\dR(X)_0=H^2_\dR(X)$.
\end{rem}

\begin{prop}\label{VDR}
Let $S^o\subset C$ be a Zariski open set such that
$U^o:=f^{-1}(S^o)\to S^o$ is smooth and the map 
\begin{equation}\label{GMB}
\ol{\nabla}:f_*\Omega^1_{U^o/S^o}\lra \Omega^1_{S^o}\ot R^1f_*\O_{U^o}
\end{equation}
induced from the Gauss-Manin connection
is bijective. Assume $S^o\ne\emptyset$ 
$($this is true if $f$ has a totally degenerate semistable fiber 
by Lem. \ref{monodromy-lem}$)$. Then the following hold.
\begin{enumerate}
\item[{\rm (1)}]
Let $S\subset C$ be an arbitrary Zariski open set and put $V=f^{-1}(S)$. 
Then $H^2_\dR(V)_0\op\NF_\dR(V)=H^2_\dR(V)$.
If $V\ne X$, then we also have $H^2_\dR(V)_0=\Image[\vg(V,\Omega^2_V)\lra H^2_\dR(V)]$.
\item[{\rm (2)}] Let $S_1\supset S_2$ and $V_i=f^{-1}(S_i)$. 
Then there is an exact sequence
\[
0\lra H^2_\dR(V_1)_0\lra H^2_\dR(V_2)_0\lra \bigoplus_{s\in S_1-S_2}H_1^\dR(f^{-1}(s)).
\]
\end{enumerate}
\end{prop}
\begin{pf}
We may assume $K=\ol{K}$ by Rem. \ref{nf-rem}.
We first prove (1). In case $V=X$, this follows from Remark \ref{vdr-nf-rem}.
Assume $V\ne X$.
We consider a spectral sequence
\[
E^{pq}_1=H^q(V,\Omega^p_{V})\Longrightarrow
H^{p+q}_\dR(V).
\]
Since $S$ is affine by the assumption, $E_1^{pq}=H^q(V,\Omega^p_{V})
=\vg(S,R^qf_*\Omega^p_{V})=0$ unless $p\leq 2$ and $q\leq 1$, so that
we have
\[
E^{20}_3=E^{20}_\infty=\Image\vg(V,\Omega^2_{V}),\quad E^{11}_2=E^{11}_\infty,
\quad E_2^{02}=0,
\]
\begin{equation}\label{VDR-0}
0\lra \Image\vg(V,\Omega^2_{V})\lra H^2_\dR(V)\lra E^{11}_\infty\lra0. 
\end{equation}
\begin{lem}\label{VDR-lem1}
The composition of maps
\[
\NF_\dR(V)\lra H^2_\dR(V)\lra E^{11}_\infty
\]
is surjective.
\end{lem}
\begin{pf}
Let $Q^o:=S\cap S^o$ and $j:V^o:=f^{-1}(Q^o)\hra V$ be the open immersion.
Consider a commutative diagram
\[
\xymatrix{
\vg(V,j_*\Omega^1_{V^o}/\Omega^1_{V})\ar[r]^{\quad \delta}&
H^1(V,\Omega^1_{V})\ar[r]^{j^*}\ar[d]^{d}&
H^1(V^o,\Omega^1_{V^o})\ar[d]^{d^o}&\mbox{(exact)}\\&
H^1(V,\Omega^2_{V})\ar[r]&
H^1(V^o,\Omega^2_{V^o})
}
\]
Let $x\in \ker d$. 
Then $j^*(x)\in \ker d^o$.
We first see that the kernel of $d^o$ is one-dimensional, generated by the cycle class 
$[e(C)]$.
Indeed, since $\ol{\nabla}$ \eqref{GMB} is bijective, 
one has $R^1f_*\Omega^1_{V^o}\os{\cong}{\lra}R^1f_*\Omega^1_{V^o/Q^o}$, 
and this is generated by the cycle class of $e(C)$ as $\O_{Q^o}$-module.
Then one can identify the map $d^o$ with $d\ot\id:\O(Q^o)\ot[e(C)]\to
\vg(Q^o,R^1f_*\Omega^2_{V^o})\cong
\vg(Q^o,\Omega^1_{Q^o}\ot R^1f_*\Omega^1_{V^o/Q^o})
=\vg(Q^o,\Omega^1_{Q^o})\ot [e(C)]$. Since the characteristic of $K$ is zero, the kernel of it
is one-dimensional over $K$. This means $\ker~d^o$ is generated by the cycle class $[e(C)]$.
Thus $x':=x-c[e(C)]$ for some $c\in K$
is contained in $\ker(j^*)=\Image~\delta$.
However, as is well-known, the image of $\delta$
is generated by the cycle classes of the irreducible components of $V-V^o$.
This shows that $x$ is a linear combination of the cycle classes 
of $e(C)$ and fibral divisors. Since $\ker(d)\to E^{11}_2=E^{11}_\infty$ is surjective,
we are done.
\end{pf}
We turn to the proof of (1).
The composition of maps
\begin{equation}\label{VDR-1}
\NF_\dR(V)\lra H^2_\dR(V) \lra \prod_{s\in S}H^2_\dR(f^{-1}(s))
\end{equation}
is given by intersection pairing, and hence is injective
by Zariski's lemma (\cite{barth} III (8.2)).
Moreover since the composition
\[
\vg(V,\Omega^2_{V})\lra H^2_\dR(V) \lra \prod_{s\in S}H^2_\dR(f^{-1}(s))
\]
is obviously zero, the second arrow in \eqref{VDR-1}
factors through $E_2^{11}=E_\infty^{11}$ (cf. \eqref{VDR-0}).
Summing up this and Lem. \ref{VDR-lem1}, we have a commutative diagram 
\begin{equation}\label{VDR-2}
\xymatrix{
0\ar[r]&\Image\vg(V,\Omega^2_{V})
\ar[r]&H^2_\dR(V/K)\ar[r]& E^{11}_\infty\ar[d]\ar[r]&0&\mbox{(exact)}\\
&&\NF_\dR(V)\ar[r]^{\subset\qquad}\ar[u]^\cup\ar[ru]^\cong&
\prod_{s\in S}  H^2_\dR(f^{-1}(s))
}
\end{equation}
with an exact row.
This shows (1).

\medskip

Next we show (2).
Let $\langle f^{-1}(s)\rangle_{s\in S_1-S_2}$
denotes the $K$-submodule of $H^2_\dR(V_1)$ generated by the cycle classes
of components of $f^{-1}(s)$ for $s\in S_1-S_2$.
By (1) we have a commutative diagram
\[
\xymatrix{
&&0\ar[d]&0\ar[d]\\
&&\NF_\dR(V_1)\ar[d]\ar[r]^{a}&\NF_\dR(V_2)\ar[d]\\
0\ar[r]&\langle f^{-1}(s)\rangle_{s\in S_1-S_2}\ar[r]&H^2_\dR(V_1)\ar[r]\ar[d]
&H^2_\dR(V_2)\ar[r]\ar[d]
&\bigoplus_{s\in S_1-S_2}H_1^\dR(f^{-1}(s))\\
&&H^2_\dR(V_1)_0\ar[d]\ar[r]&H^2_\dR(V_2)_0\ar[d]\\
&&0&0
}
\]
with exact row and columns.
The map $a$ is surjective and $\ker(a)$ is onto 
$\langle f^{-1}(s)\rangle_{s\in S_1-S_2}$.
Now the desired assertion follows from the snake lemma.
\end{pf}

%%%%%%%%%%%%%%%%%%%%%%%%%%%%%%%%%%%%%%%%%%%%%%%

\subsection{Deligne's canonical extension}\label{hodge-sect}
Let $j:S\hra C$ be a Zariski open set such that $U=f^{-1}(S)\to S$ is smooth.
Put $T:=C-S$ and $D:=f^{-1}(T)$.
By taking the embedded resolution of singularities if necessary, 
we can assume that $D_{\mathrm{red}}$ is a NCD.
We then consider the de Rham cohomology groups
\[
H^q_\dR(U)=H^q_\zar(U,\Omega^\bullet_{U})
\cong
H^q_\zar(X,\Omega^\bullet_X(\log D))
\]
with the Hodge filtration 
\[
F^p H^q_\dR(U):=\Image[H^q(X,\Omega^{\bullet\geq p}_X(\log D))\hra
H^q(X,\Omega^\bullet_X(\log D))].
\]
Define a sheaf $\Omega^1_{X/C}(\log D)$ by the exact sequence
\[
0\lra 
f^*\Omega^1_{C}(\log T)\lra
\Omega^1_{X}(\log D)\lra
\Omega^1_{X/C}(\log D)\lra 0.
\]
This is a locally free sheaf of rank one.
Let
$\cH_e:=R^1f_*\Omega^\bullet_{X/C}(\log D)$
be {\it Deligne's canonical extension} and 
\begin{equation}\label{d-hf0}
\cH_e^{1,0}:=f_*\Omega^1_{X/C}(\log D),\quad
\cH_e^{0,1}:=\cH_e/\cH_e^{1,0}\cong R^1f_*\O_X
\end{equation}
the Hodge filtration (cf. Appendix \S \ref{de-sect}).
The {\it Gauss-Manin connection}
\begin{equation}\label{gauss-manin}
\nabla:\cH_e\lra \Omega^1_C(\log T)\ot\cH_e
\end{equation}
is defined to be the connecting homomorphism arising from an exact sequence  
\begin{equation}\label{hodge-exmp-12}
0\to
f^*\Omega^1_{C}(\log T)\ot \Omega^{\bullet-1}_{X/C}(\log D)\to
\Omega^\bullet_{X}(\log D)\to
\Omega^\bullet_{X/C}(\log D)\to 0
\end{equation}
(see Appendix \eqref{gm-diag} for a remark on sign.)
Note that $\cH_e$ is characterized as a subbundle of $j_*\cH$ 
such that the eigenvalues of $\Res(\nabla)$ are in $[0,1)$.
In particular it does not depend on the choice of $D$.
Write
\[
H^q_\dR(C,\cH_e):=H^q_\zar(C,\cH_e\to \Omega^1_C(\log T)\ot\cH_e).
\]
\begin{thm}[cf. \cite{SZ} \S 5]\label{Hodge}
There is the natural isomorphism
\begin{equation}\label{hodge-exmp-8}
H^1_\dR(C,\cH_e)\os{\cong}{\lra}
H^2_\dR(U)_0.
\end{equation}
Moreover 
the Hodge filtration corresponds in the following way. 
\begin{align}
&F^1H^2_\dR(U)_0\cong 
H^1_\zar(C,\cH^{1,0}_e
\to \Omega^1_C(\log T)\ot\cH^{1,0}_e)
\rangle\label{hodge-exmp-9}\\
&F^2H^2_\dR(U)_0\cong H^0_\zar(C,\Omega^1_C(\log T)\ot\cH^{1,0}_e)\label{hodge-exmp-10}\\
&\mathrm{Gr}_F^0H^2_\dR(U)_0\cong H^1_\zar(C,\cH^{0,1}_e)\label{hodge-exmp-11}
\end{align}
\end{thm}
\begin{pf}
The exact sequence
\eqref{hodge-exmp-12}
gives rise to a spectral sequence
\[
E_2^{pq}=H^p_\dR(C,R^qf_*\Omega^{\bullet}_{X/C}(\log D))\Longrightarrow
H^{p+q}_\dR(U).
\]
This yields 
\[
0\lra H^1_\dR(C,\cH_e)
\lra H^2_\dR(U)\lra H^0_\dR(C,R^2f_*\Omega^{\bullet}_{X/C}(\log D)) 
\lra 0.
\]
Since the last term is one-dimensional, isomorphic to $H^1_\dR(X_t)$, 
we have \eqref{hodge-exmp-8}.

\eqref{hodge-exmp-12} induces an exact sequence 
\[
0\to f^*\Omega^1_C(\log T)\ot
\Omega^{\bullet-1\geq p-1}_{X/C}(\log D)
\to 
\Omega^{\bullet\geq p}_X(\log D)\lra 
\Omega^{\bullet\geq p}_{X/C}(\log D)
\to 0
\]
and this yields
\[
\xymatrix{
H^1_\zar(C,R^1f_*\omega_{X/C}^{\bullet\geq p}
\to \Omega^1_C(\log T)\ot R^1f_*\omega^{\bullet\geq p-1}_{X/C})
\ar[d]\\
H^2_\zar(X,\omega^{\bullet\geq p}_{X/C}
\to f^*\Omega^1_C(\log T)\ot\omega^{\bullet\geq p-1}_{X/C})
\ar[r]^{\qquad\quad\cong}& H^2_\zar(X,\Omega^{\bullet\geq p}_X(\log D))\\
}
\]
where $\omega_{X/C}^\bullet:=\Omega^\bullet_{X/C}(\log D)$.
Now \eqref{hodge-exmp-9}, 
\eqref{hodge-exmp-10} and \eqref{hodge-exmp-11} easily follow from this.
\end{pf}
\begin{lem}\label{monodromy-lem}
Let 
\begin{equation}\label{PF-rem-eq}
\ol{\nabla}:\cH_e^{1,0}\lra \Omega^1_{C}(\log T)\ot \cH_e^{0,1}
\end{equation}
be the $\O_C$-linear map
induced from the Gauss-Manin connection \eqref{gauss-manin}.
Let $f^{-1}(P)$ be a fiber over $P\in T$.
If $f^{-1}(P)$ is semistable, 
then $\ol{\nabla}$ is bijective on a neighborhood of $P$
if and only if $f^{-1}(P)$ is totally degenerate.
If there is a component of $f^{-1}(P)$ which is not rational, then 
$\ol{\nabla}$ is not surjective on a neighborhood of $P$.
\end{lem}
\begin{pf}
Since the both sides of \eqref{PF-rem-eq} are locally free sheaves of the same rank,
the bijectivity of \eqref{PF-rem-eq} is equivalent to the surjectivity of it.
By Nakayama's lemma, it is also equivalent 
to the surjectivity modulo the maximal ideal at $P$. 
We may assume $K=\C$. Let $t\in \O_{C,P}$ be the uniformizer and write $\C_P=\O_{C,P}/t
\O_{C,P}$.
There is an isomorphism
$\cH_e\ot\C_p\cong H^2_\dR(X_t)$ 
where $X_t$ is a smooth fiber 
(which is ``close to $f^{-1}(P)$") (\cite{steenbrink} (2.18)).
Moreover $\Res_P(\nabla):H^2_\dR(X_t)\to H^2_\dR(X_t)$ coincides with
the log monodromy operator $N$ such that the eigenvalues of $N$ are in $[0,1)$
(\cite{steenbrink} (2.21)).
Let $\hat{F}^\bullet$ be the filtration on $H^2_\dR(X_t)$ induced from \eqref{d-hf0}.
Then $\ol{\nabla}$ \eqref{PF-rem-eq} is surjective on a neighborhood of $P$
if and only if the map
\[
\ol{N}:\hat{F}^1 H^1_\dR(X_t)\lra \Gr^0_{\hat{F}}H^1_\dR(X_t)
\]
induced from $N$ is surjective, or equivalently injective.

Assume that $f^{-1}(P)$ is a semistable fiber.
Then $\hat{F}^\bullet$ is the limiting Hodge filtration due to Steenbrink \cite{steenbrink}.
One has $\ker(\ol{N})=F^1H^1_\dR(f^{-1}(P))$
by the local invariant cycle theorem (\cite{steenbrink} (5.12)).
Therefore $\ker(\ol{N})=0$ if and only if
$f^{-1}(P)$ is a totally degenerate curve.

There is an obvious inclusion
$F^1H^1_\dR(f^{-1}(P))\hra \ker(N)\cap \hat{F}^1\subset \ker(\ol{N})$
without the assumption that $f^{-1}(P)$ is a semistable fiber.
Therefore if $f^{-1}(P)$ contains a non-rational curve, then
$\ker(\ol{N})\ne0$.
\end{pf}
\begin{rem}
In case of elliptic fibration, it follows from Thm. \ref{A1} and \ref{A2}
that
$\ol{\nabla}$ \eqref{PF-rem-eq} is bijective if and only if either of the following 
conditions holds.
\begin{enumerate}
\item[(i)]
$f^{-1}(P)$ is a (non-smooth) semistable fiber (i.e. multiplicative),
\item[(ii)]
$f^{-1}(P)$ is additive and 
\[
\Res_P\left(\frac{t(2g_2dg_3-3g_3dg_2)}{\Delta}\right) \ne 0,\quad \Delta:=g_2^3-27g_3^2
\]
where $t\in \O_{C,P}$ is a uniformizer and
$y^2=4x^3-g_2x-g_3$ is the minimal Weierstrass equation of $f$ over a neighborhood of $P$.
\end{enumerate}   
However, it seems difficult to give a complete criterion of
the bijectivity of $\ol{\nabla}$ in case $g>1$.
\end{rem}

%%%%%%%%%%%%%%%%%%%%%%%%%%%%%%%%%%%%%%%%%%%%%%%%%
\subsection{Relative cohomology}\label{relative-sect}
For a smooth manifold $M$, we denote by 
$\cA^q(M)$ the space of smooth differential $q$-forms on $M$
with coefficients in $\C$.

\medskip

Let $f:X\to C$ be a fibration of curves over $\C$. 
Let $S\subset C$ be an arbitrary Zariski open set, and put $V:=f^{-1}(S)$.
Let $D\subset V$ be a fiber.
Let $\rho:\wt{D}\to D$ be the normalization and $\Sigma\subset D$
the set of singular points. Let $s:\wt{\Sigma}:=\rho^{-1}(\Sigma)
\hra\wt{D}$ be the inclusion.
There is the exact sequence
\[
0\lra \O_{D}\os{\rho^*}{\lra} \O_{\wt{D}}\os{s^*}{\lra} 
\C_{\wt{\Sigma}}/\C_\Sigma\lra 0
\]
where $\C_{\wt{\Sigma}}=\mathrm{Maps}(\wt{\Sigma},\C)=\Hom(\Z\wt{\Sigma},\C)$, 
$\rho^*$ and $s^*$ are the pull-back.
We define $\cA^\bullet(D)$ to be the mapping fiber of 
$s^*:\cA^\bullet(\wt{D})\to \C_{\wt{\Sigma}}/\C_\Sigma$:
\[
\cA^0(\wt{D})
\os{s^*\op d}{\lra} \C_{\wt{\Sigma}}/\C_\Sigma\op \cA^1(\wt{D})
\os{0\op d}{\lra} \cA^2(\wt{D})
\]
where the first term is placed in degree 0.
Then
\[
H^q_\dR(D)=H^q(\cA^\bullet(D))
\]
is the de Rham cohomology of $D$, which fits into the exact sequence
\[
\cdots\lra H^0_\dR(\wt{D})\lra \C_{\wt{\Sigma}}/\C_\Sigma
\lra H^1_\dR(D)\lra H^1_\dR(\wt{D})\lra\cdots.
\]
There is the natural pairing 
\begin{equation}\label{pairingD}
H_1(D,\Z)\otimes H^1_\dR(D)\lra \C,
\quad
\gamma\ot z\mapsto\int_\gamma z:=\int_\gamma\eta-c(\partial(\rho^{-1}\gamma))
\end{equation}
where $z=(c,\eta)\in \C_{\wt{\Sigma}}/\C_\Sigma\op\cA^1(\wt{D})$ with $d\eta=0$
and $\partial$ denotes the boundary of homology cycles.

\medskip

We define $\cA^\bullet(V,D)$
to be the mapping fiber of 
$j^*:\cA^\bullet(V)\to \cA^\bullet(D)$ the pull-back of $j:D\hra V$:
\[
\cA^0(V)\os{\cD_0}{\lra} \cA^0(\wt{D})\op \cA^1(V)
\os{\cD_1}{\lra} \C_{\wt{\Sigma}}/\C_\Sigma\op\cA^1(\wt{D})\op \cA^2(V)\os{\cD_2}{\lra}\cdots
\]
where 
\[
\cD_0=j^*\op d,\quad
\cD_1=\begin{pmatrix}
-(s^*\op d)&j^*\\
&d
\end{pmatrix},\quad
\cD_2=\begin{pmatrix}
-(0\op d)&j^*\\
&d
\end{pmatrix},\ldots
\]
Then 
\begin{equation}\label{exp-9-0}
H^q_\dR(V,D)=H^q(\cA^\bullet(V,D))
\end{equation}
is the de Rham cohomology which
fits into the exact sequence
\begin{equation}\label{exp-9}
\cdots\lra H^{q-1}_\dR(D)\lra H^q_\dR(V,D)\lra H^q_\dR(V)\lra H^q_\dR(D)\lra\cdots. 
\end{equation}
An element of $H^2_\dR(V,D)$ is described by $z=(c,\eta,\omega)
\in \C_{\wt{\Sigma}}/\C_\Sigma\op\cA^1(\wt{D})\op \cA^2(V)$ with
$j^*\omega=d\eta$ and $d\omega=0$ which are subject to
relations $(s^*f, df,0)=0$ and $(0,j^*\theta,d\theta)=0$ for $f\in \cA^0(\wt{D}_0)$
and $\theta\in \cA^1(V)$.
The natural pairing
\begin{equation}\label{pairingVD1}
H_2(V,D;\Z)\ot H_\dR^2(V,D)\lra \C,\quad \Gamma\ot z
\longmapsto \int_\Gamma z
\end{equation}
is given by
\begin{equation}\label{pairingVD2}
\int_\Gamma z:=\int_\Gamma \omega-\int_{\partial\Gamma}(c,\eta)
=\int_\Gamma \omega-\int_{\partial\Gamma}\eta+c(\rho^{-1}(\partial\Gamma)).
\end{equation}

\begin{lem}\label{pairingVD4}
Put
\begin{equation}\label{pairingVD4-1}
H^2_\dR(V,D)_0:=\ker[H^2_\dR(V,D)\lra \prod_{s\in S}H^2_\dR(f^{-1}(s))\times H^2_\dR(e(S))],
\end{equation}
and hence there is an exact sequence
\begin{equation}\label{pairingVD4-2}
F^1H^1_\dR(D)\lra F^1H^2_\dR(V,D)_0\lra F^1H^2_\dR(V)_0\lra0.
\end{equation}
For $\omega\in F^1H^2_\dR(V)_0$ let ${\omega}_{V,D}\in F^1H^2_\dR(V,D)_0$ be a lifting.
Let $\Gamma\in H_2(V,D;\Q)$.
If $\gamma:=\partial\Gamma\in H_1^B(D,\Q)$ belongs to the Hodge $(0,0)$-part, then
$\int_\Gamma{\omega}_{V,D}$ does not depend on the choice of the lifting ${\omega}_{V,D}$.
\end{lem}
\begin{pf}
Let ${\omega}_{V,D}^\prime$ be another lifting, then 
${\omega}_{V,D}^\prime-{\omega}_{V,D}$ belongs to $F^1H_\dR(D)$.
Therefore by \eqref{pairingVD2} the assertion follows from
the fact that the pairing
\[
F^0H_1(D,\C)\ot F^1H^1_\dR(D)\to \C\] 
induced from \eqref{pairingD} is zero.
\end{pf}

\def\goodua{\Lambda^1({U})_\mathrm{rat}}
\def\goodub{\Lambda^2(U)_\mathrm{rat}}
\def\goodxa{\Lambda^1(X)_\mathrm{rat}}
\def\goodxb{\Lambda^2(X)_\mathrm{rat}}
\subsection{$\Lambda(U)_\mathrm{rat}$ and $\Lambda(X)_\mathrm{rat}$}\label{good-sect}
Let $j:S\hra C$ be a Zariski open set such that $U=f^{-1}(S)\to S$ is smooth.
Put $T:=C-S$ and $D:=f^{-1}(T)$. We denote by $\cH_e$ Deligne's canonical extension
as in \S \ref{hodge-sect}.
Let $C^o\subset C$ be the maximal open set such that
\eqref{PF-rem-eq} is bijective on $C^o$.
We assume $C^o\ne\emptyset$.
By Lemma \ref{monodromy-lem}, if $f$ has a totally degenerate semistable fiber,
then $P\in C^o$ and hence $C^o\ne\emptyset$. 
Put $X^o:=f^{-1}(C^o)$.
We first introduce two spaces of rational 2-forms
\[
\goodub\subset\goodua\subset\vg(C^o,\Omega^1_{C}(\log T)\ot\cH_e^{1,0})
\subset\vg(U\cap X^o,\Omega^2_X).
\]
Define
\begin{align*}
\goodub&:=\Image[\vg(C,\Omega^1_{C}(\log T)\ot\cH_e^{1,0})
\hra \vg(C^o,\Omega^1_{C}(\log T)\ot\cH_e^{1,0})]\\
&\cong F^2H^2_\dR(U).
\end{align*}
We define $\goodua$ in the following way.
Let us consider a diagram
\begin{equation}\label{hodge-pf1}
\xymatrix{
&0\ar[d]\\
&\Omega^1_{C}(\log T)\ot\cH_e^{1,0}|_{C^o}\ar[d]\\
\cH_e^{1,0}|_{C^o}\ar[r]^{\nabla\qquad\quad}\ar[d]_=&\Omega^1_C(\log T)\ot\cH_e|_{C^o}\ar[d]\\
\cH_e^{1,0}|_{C^o}\ar[r]^{\ol{\nabla}\qquad\quad}&\Omega^1_C(\log T)\ot\cH_e^{0,1}|_{C^o}\ar[d]\\
&0}
\end{equation}
Since $\ol{\nabla}$ is bijective by definition of $C^o$,
one has an isomorphism
\begin{equation}\label{hodge-pf2}
H^1_\zar(C^o,\cH_e^{1,0}\to\Omega^1_C(\log T)\ot\cH_e)
\os{\cong}{\longleftarrow}
\vg(C^o,\Omega^1_C(\log T)\ot\cH_e^{1,0}).
\end{equation}
We define $\goodua$ to be the image of the composition of the following maps
\begin{align}
H^1_\zar(C,\cH_e^{1,0}\to\Omega^1_C(\log T)\ot\cH_e)
&\lra
H^1_\zar(C^o,\cH_e^{1,0}\to\Omega^1_C(\log T)\ot\cH_e)\label{hodge-pf2-0}\\
&\us{\cong}{\os{\eqref{hodge-pf2}}{\longleftarrow}}
\vg(C^o,\Omega^1_C(\log T)\ot\cH_e^{1,0}).\label{hodge-pf2-1}
\end{align}
\begin{prop}\label{good-prop1}
\[
F^1H^2_\dR(U)_0\cong
H^1_\zar(C,\cH_e^{1,0}\to\Omega^1_C(\log T)\ot\cH_e)
\os{\cong}{\lra}\goodua.
\]
\end{prop}
\begin{pf}
The first isomorphism is due to Thm.\ref{Hodge}.
To show the second, it is enough to show the injectivity of
\eqref{hodge-pf2-0}. 
However this follows from the fact that $F^1H^2_\dR(U)_0\to F^1H^2(U\cap X^o)_0$ is 
injective by Prop. \ref{VDR} (2).
\end{pf}
Let
$\Res_D:H^2_\dR(U)\lra H_1^\dR(D)$
be the residue map along $D$.
We define
\begin{equation}\label{hodge-pf5}
\Lambda^i(X)_\mathrm{rat}:=\Lambda^i(U)_\mathrm{rat}\cap\ker(\Res_D).
\end{equation} 
By definition one has 
\begin{equation}\label{hodge-pf6}
\Lambda^i(X)_\mathrm{rat}\subset \vg(C^o,\Omega^1_C(\log T)\ot\cH_e^{1,0})\cap\ker(\Res_D)
=\vg(X^o,\Omega^2_{X^o}).
\end{equation}
Moreover $\Lambda^i(X)_\mathrm{rat}$ does not depend on the choice of $U$.
\begin{prop}\label{good-prop2}
$\Lambda^i(X)_\mathrm{rat} \cong F^iH^2_\dR(X)_0.$
In particular, it is stable under a birational transformation $X'\to X$.
\end{prop}
\begin{pf}
Prop. \ref{VDR} (2) and 
the definition of $\Lambda^i(X)_\mathrm{rat}$ give rise to
a commutative diagram
\[
\xymatrix{
0\ar[r]&\Lambda^i(X)_\mathrm{rat}\ar[d]\ar[r]
&\Lambda^i(U)_\mathrm{rat}\ar[d] \ar[r] &H_{1,\dR}(D)\ar@{=}[d]\\
0\ar[r]&F^iH^2_\dR (X)_0 \ar[r]&F^iH^2_\dR (U)_0\ar[r]&H_{1,\dR}(D)
}
\]
with exact rows. Now the assertion follows from 
Prop. \ref{good-prop1}.
\end{pf}
The following theorem is one of the technical key results,
will be used in the proof of the main theorem (see Thm. \ref{ExtMHS} (2)).
\begin{thm}\label{hodge-mot}
Let $D_0=f^{-1}(P)$ be a fiber contained in $X^o$.
For $\omega\in \goodxa$, let $\hat{\omega}=(0,0,\omega)\in H^2_\dR(X^o,D_0)$ be the 
cohomology class in terms of the RHS of \eqref{exp-9-0}.
Then $\hat{\omega}$ 
belongs to $F^1H^2_\dR(X^o,D_0)$
where $F^\bullet$ denotes the Hodge filtration.
\end{thm}

To prove Theorem \ref{hodge-mot}, we may assume that $(D_0)_\mathrm{red}$ 
and $D_\mathrm{red}$ are NCD's.
Let $j:D_0\hra X$ and $\rho:\wt{D}_0\to D_0$ the normalization.
Let $\Sigma\subset D_0$ be the singular locus and put 
$s:\wt{\Sigma}:=\rho^{-1}(\Sigma)\hra \wt{D}_0$.

We use the Cech cocycles to describe the de Rham cohomology groups.
Let us denote by $(\check{C}^\bullet({\mathscr F}),\delta)=
(\check{C}^\bullet(X,{\mathscr F}),\delta)$ 
the Cech complex of a Zariski sheaf $\mathscr F$.
Let
\[
\check{C}^1(\O_X)\times
\check{C}^0(\Omega^1_X)
\os{\cD_1}{\lra}
\check{C}^2(\O_X)\times
\check{C}^1(\Omega^1_X)
\times
\check{C}^0(\Omega^2_X)
\os{\cD_2}{\lra}
\check{C}^3(\O_X)\times
\check{C}^2(\Omega^1_X)
\times
\check{C}^1(\Omega^2_X)
\]
\[
\cD_1=
\begin{pmatrix}
\delta&-d&\\
&\delta&d
\end{pmatrix}
,\quad
\cD_2=
\begin{pmatrix}
\delta&d&\\
&\delta&-d\\
&&\delta
\end{pmatrix}.
\]
Then the cohomology of the middle term of the above complex gives $H^2_\dR(X)$.
In the same way we obtain the description of $H^2_\dR(U)=H^2(X,\Omega^\bullet_X(\log D))$.
Let
\begin{multline*}
\check{C}^1(\O_V)\times
\check{C}^0(\O_{\wt{D}_0}\op\Omega^1_V)
\os{\cD_3}{\lra}
\check{C}^2(\O_V)\times
\check{C}^1(\O_{\wt{D}_0}\op\Omega^1_V)
\times
\check{C}^0(\O_{\wt{\Sigma}}/\O_\Sigma\op\Omega^1_{\wt{D}_0}\op\Omega^2_V)\\
\os{\cD_4}{\lra}
\check{C}^3(\O_V)\times
\check{C}^2(\O_{\wt{D}_0}\op\Omega^1_V)
\times
\check{C}^1(\O_{\wt{\Sigma}}/\O_\Sigma\op\Omega^1_{\wt{D}_0}\op\Omega^2_V)
\end{multline*}
\[
\cD_3=
\begin{pmatrix}
\delta&-(j^*\op d)&\\
&\delta&T
\end{pmatrix},\quad
\cD_4=
\begin{pmatrix}
\delta&j^*\op d&\\
&\delta&-T\\
&&\delta
\end{pmatrix},\quad
T=
\begin{pmatrix}
-s^*&-d&\\
&j^*&d
\end{pmatrix}.
\]
Then the cohomology of the middle term of the above complex gives $H^2_\dR(X,D_0)$.
Note that 
\[
\hat{\omega}=(0)\times(0,0)\times(0,0,\omega)\in H^2_\dR(X^o,D_0)
\]
in terms of the Cech cocycles.
\begin{lem}\label{extra2-1}
Let $z=(0)\times (\alpha_{ij})\times(\beta_i)\in 
\ker(\cD_2)$
be a Cech cocycle, and $[z]\in F^1H^2_\dR(X)$
its cohomology class $[z]\in F^1H^2_\dR(X)$.
Assume $[z]\in \ker[H^2_\dR(X)\to H^2_\dR(D_0)\cong H^2_\dR(\wt{D}_0)]$, so that
there is $(\epsilon_i)\in \check{C}^0(\Omega^1_{\wt{D}_0})$ such that
$\alpha_{ij}|_{\wt{D}_0}=\epsilon_j-\epsilon_i$.
Put
\[
z_{X,D_0}:=(0)\times (0,\alpha_{ij})\times(0,\epsilon_i,\beta_i)\in \ker(\cD_4).
\]
Then we have $[z_{X,D_0}]\in F^1H^2_\dR(X,D_0)$ and this is a lifting of $[z]$ via the map
$F^1H^2_\dR(X,D_0)\to F^1H^2_\dR(X)$.
\end{lem}
\begin{pf}
Obvious from the definition of Hodge filtration.
\end{pf}

We turn to the proof of Theorem \ref{hodge-mot}\footnote
{One cannot directly apply Lemma \ref{extra2-1} to a Cech cocycle
$(0)\times (0)\times (\omega)\in 
\check{C}^2(\O_{X^o})\times
\check{C}^1(\Omega^1_{X^o})
\times
\check{C}^0(\Omega^2_{X^o})$ to show $\hat{\omega}\in F^1H^2_\dR(X^o,D_0)$ 
because $X^o$ is not complete.}. 
Let 
\[
(\alpha_{ij})\times (\beta_i)\in \check{C}^1(\cH^{1,0}_e)\times
\check{C}^0(\Omega^1_C(\log T)\ot\cH_e)
\]
be a corresponding Cech cocycle to $\omega\in \goodxa$, and this defines
\[
z:=(0)\times (\eta_{ij})\times (\pi_i)\in \check{C}^2(\O_X)\times
\check{C}^1(\Omega^1_X(\log D))\times
\check{C}^0(\Omega^2_X(\log D))
\]
in a natural way.
Since $[z]\in F^1H^2_\dR(U)$ lies in the image of $F^1H^2_\dR(X)_0$,
there is a Cech cocycle
$
w=(0)\times(*)\times(*)\in \check{C}^2(\O_X)\times
\check{C}^1(\Omega^1_X)\times
\check{C}^0(\Omega^2_X)
$
such that $[w]\in F^1H^2_\dR(X)_0$ and $[w]|_U=[z]$ in $H^2_\dR(U)$.
Since the map 
$H^2(X,\Omega^{\bullet\geq 1}_X(\log D))\to
H^2(X,\Omega^{\bullet}_X(\log D))$ is injective, we see that 
there is $\wt{y}=(0)\times (\wt{\nu}_i)\in \check{C}^1(\O_X)\times
\check{C}^0(\Omega^1_X(\log D))$
such that 
\begin{equation}\label{hodge-pf9}
w=z-\cD_1(\wt{y})=
(0)\times (\eta_{ij}-(\wt{\nu}_j-\wt{\nu}_i))\times (\pi_i-d\wt{\nu}_i)
\end{equation}
and this belongs to $\check{C}^2(\O_X)\times
\check{C}^1(\Omega^1_X)\times
\check{C}^0(\Omega^2_X)$.
Therefore, by Lemma \ref{extra2-1}
there is $(\epsilon_i)\in \check{C}^0(\Omega^1_{\wt{D}_0})$ such that
$\epsilon_j-\epsilon_i=\eta_{ij}-(\wt{\nu}_j-\wt{\nu}_i)|_{\wt{D}_0}$, and
\begin{equation}\label{hodge-pf9-0}
z_{X,D_0}:=
(0)\times (0,\eta_{ij}-(\wt{\nu}_j-\wt{\nu}_i))\times (0,\epsilon_i,\pi_i-d\wt{\nu}_i)
\end{equation}
defines a lifting $[z_{X,D_0}]\in F^1H^2_\dR(X,D_0)$ of $\omega\in \goodxa$. 
We want to show that
$z_{X,D_0}\equiv \hat{\omega}=(0)\times (0,0)\times(0,0,\omega)$ in $H^2_\dR(X,D_0)$
modulo the image of $F^1H^1_\dR(D_0)$.
To do this it is enough to show 
\begin{equation}\label{hodge-pf9-2}
z_{X,D_0}|_V= \hat{\omega}|_V\in H^2_\dR(V,D_0)/\Image F^1H^1_\dR(D_0)
\end{equation}
for a (sufficiently small) neighborhood $V$ of $D_0$.

By the definition of the locus $C^o$, there is 
$y_0\in \check{C}^0(\cH_e|_{C^o})$ such that 
\[
(0)\times (\omega)=(\alpha_{ij})\times (\beta_i)-\cD_0(y_0)
\]
where $\cD_0:\check{C}^0(\cH_e)\to
\check{C}^1(\cH^{1,0}_e)\times
\check{C}^0(\Omega^1_C(\log T)\ot\cH_e)$.
This means that there is $y=(0)\times(\nu_i)
\in \check{C}^1(\O_{X^o})\times\check{C}^0(\Omega^1_{X^o}(\log D))$
such that
\begin{equation}\label{hodge-pf8}
z|_{X^o}-\cD_1(y)=
(0)\times (\eta_{ij}|_{X^o}-(\nu_j-\nu_i))\times (\pi_i|_{X^o}-d\nu_i)
=(0)\times(0)\times (\omega).
\end{equation}
Therefore we have
\begin{equation}\label{hodge-pf9-1}
z_{X,D_0}|_{X^o}=
(0)\times (0,(\nu_j-\wt{\nu}_j)-(\nu_i-\wt{\nu}_i))\times (0,\epsilon_i,
\omega+d\nu_i-d\wt{\nu}_i)
\end{equation}
with
\[
((\nu_j-\wt{\nu}_j)-(\nu_i-\wt{\nu}_i))|_{\wt{D}_0}=\epsilon_j-\epsilon_i.
\]
We note that $\nu_i$ and $\wt{\nu}_i$ have at most log pole along $D_0$.
\begin{lem}\label{extra2}
Let $V$ be a (sufficiently small) neighborhood $D_0$.
Let $t\in \O_{C,P}$ be a uniformizer at $P$.
Then there is a constant $c$ such that
\[
\theta_i:=\wt{\nu}_i|_V-\nu_i|_V-c\frac{dt}{t}
\]
has no log pole along $D_0$.
\end{lem}
\begin{pf}
There is the exact sequence
\[
0\lra\Omega^1_V\lra
\Omega^1_V(\log D_0)\os{\Res}{\lra} \O_{\wt{D}_0}\lra 0.
\]
Since $((\nu_j-\wt{\nu}_j)-(\nu_i-\wt{\nu}_i))$ has log pole along $D_0$, 
one has $\Res(\nu_i-\wt{\nu}_i)
=\Res(\nu_j-\wt{\nu}_j)$ and hence it defines
\[
e:=(\Res(\nu_i-\wt{\nu}_i))_i\in 
\ker[\check{C}^0(\O_{\wt{D}_0})\to \check{C}^1(\O_{\wt{D}_0})]=H^0(\O_{\wt{D}_0}).
\]
Put
$e':=((\nu_j-\wt{\nu}_j)-(\nu_i-\wt{\nu}_i))\in \check{C}^1(\Omega^1_V)$.
Then the cohomology class $[e']\in H^1(V,\Omega^1_V)$ is the image of $e$ via
the connecting homomorphism $H^0(\O_{\wt{D}_0})\to H^1(\Omega^1_V)$.
On the other hand, it follows from  \eqref{hodge-pf9-1}
that the class $[e']|_{\wt{D}_0}\in H^1(\Omega^1_{\wt{D}_0})
\cong H^2_\dR(\wt{D}_0)$ coincides with
the image of $z_{X^o,D_0}|_V$ via the composition of maps
$H^2_\dR(V,D_0)\to H^2_\dR(V)
\to H^2_\dR(\wt{D}_0)$, and hence the image of $\omega|_V$ via 
$H^2_\dR(V)\to H^2_\dR(\wt{D}_0)$. However this is obviously zero.
Thus we have
\[
e\in \ker[H^0(\O_{\wt{D}_0})\lra H^2_\dR(\wt{D}_0)]=\langle\Res(\frac{dt}{t})\rangle
\cong K
\]
where the middle equality follows from Zariski's lemma (\cite{barth} III (8.2)).
This means that there is a constant $c$ such that
\[
\theta_i:=\nu_i|_V-\wt{\nu}_i|_V-c\frac{dt}{t}
\]
has no log pole along $D_0$.
\end{pf}
Let us prove \eqref{hodge-pf9-2}.
By Lemma \ref{extra2}, one can put $\epsilon_i:=\theta_i|_{\wt{D}_0}$.
Hence we have from \eqref{hodge-pf9-1}
\begin{align*}
z_{X,D_0}|_V&=
(0)\times(0,-(\theta_j-\theta_i))
\times
(0,-\theta_i|_{\wt{D}_0},\omega|_V-d\theta_i)\\
&\equiv
(0)\times(0,0)\times(0,0,\omega|_V)\mod \Image\cD_3
\end{align*}
as required.
This completes the proof of Thm. \ref{hodge-mot}.

\subsection{Lefschetz thimbles}\label{exp-sect2}
Suppose $K=\C$. 
Let $\ol{S}\subset C$ be an arbitrary Zariski open set, and put $\ol{U}:=f^{-1}(\ol{S})$.
Put
\[
\NF^B(\ol{U}):=\Image[\bigoplus_{s\in \ol{S}}H_2(f^{-1}(s),\Z)\op
H_2(e(\ol{S}),\Z)\lra H_2(\ol{U},\Z)],
\]
\[
H_2(\ol{U},\Z)_0:=H_2(\ol{U},\Z)/\NF^B(\ol{U}).
\]
Note $H_2(e(\ol{S}),\Z)=0$ unless $\ol{S}=C$.  
By definition $H^2_\dR(\ol{U})_0\cong \Hom(H_2(\ol{U},\Q),\C)$.

\medskip

Let $S\subset \ol{S}$ be a Zariski open set such that
$U:=f^{-1}(S)\to S$ is smooth.
We put $T:=\ol{S}-S$, $D:=f^{-1}(T)$. 
Take a path $\gamma:[0,1]\to \ol{S}(\C)$, $t\mapsto\gamma_t$
such that $\gamma_t\in S(\C)$ for $t\ne0,1$.
Take a cycle $\varepsilon\in H_1(f^{-1}(\gamma_{t_0}),\Z)$ for some (fixed) $t_0\in [0,1]$.
Then it extends to a flat section 
$\varepsilon_{t}\in H_1(f^{-1}(\gamma_{t}),\Z)$ over $t\in [0,1]$
in a unique way.
Let $\Gamma(\varepsilon,\gamma)$ be the fibration over the path $\gamma$ 
whose fiber is $\varepsilon_t$. 

%WinTpicVersion4.26
\unitlength 0.1in
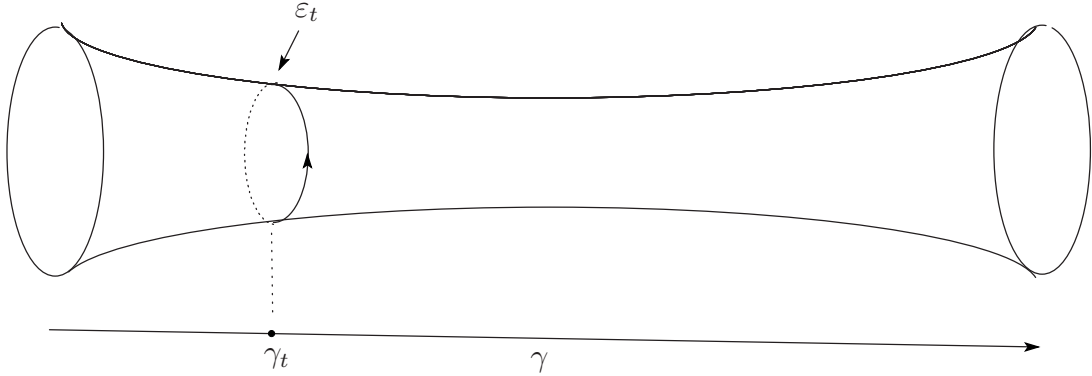
\begin{picture}( 56.1000, 19.7500)(  5.1000,-34.7500)
% ELLIPSE 2 0 3 0 Black White
% 4 3320 1740 5850 2150 780 1760 5670 1780
% 
{\color[named]{Black}{%
\special{pn 8}%
\special{ar 3320 1740 2530 410  0.1047219  3.0932053}%
}}%
% ELLIPSE 2 0 3 0 Black White
% 4 3320 1740 5850 2150 780 1760 5670 1780
% 
{\color[named]{Black}{%
\special{pn 8}%
\special{ar 3320 1740 2530 410  0.1047219  3.0932053}%
}}%
% ELLIPSE 2 0 3 0 Black White
% 4 3320 1740 5850 2150 780 1760 5670 1780
% 
{\color[named]{Black}{%
\special{pn 8}%
\special{ar 3320 1740 2530 410  0.1047219  3.0932053}%
}}%
% ELLIPSE 2 0 3 0 Black White
% 4 3320 1740 5850 2150 780 1760 5670 1780
% 
{\color[named]{Black}{%
\special{pn 8}%
\special{ar 3320 1740 2530 410  0.1047219  3.0932053}%
}}%
% ELLIPSE 2 0 3 0 Black White
% 4 3320 1740 5850 2150 780 1760 5670 1780
% 
{\color[named]{Black}{%
\special{pn 8}%
\special{ar 3320 1740 2530 410  0.1047219  3.0932053}%
}}%
% ELLIPSE 2 0 3 0 Black White
% 4 3320 1740 5850 2150 780 1760 5670 1780
% 
{\color[named]{Black}{%
\special{pn 8}%
\special{ar 3320 1740 2530 410  0.1047219  3.0932053}%
}}%
% ELLIPSE 2 0 3 0 Black White
% 4 3320 1740 5850 2150 780 1760 5670 1780
% 
{\color[named]{Black}{%
\special{pn 8}%
\special{ar 3320 1740 2530 410  0.1047219  3.0932053}%
}}%
% ELLIPSE 2 0 3 0 Black White
% 4 3320 1740 5850 2150 780 1760 5670 1780
% 
{\color[named]{Black}{%
\special{pn 8}%
\special{ar 3320 1740 2530 410  0.1047219  3.0932053}%
}}%
% ELLIPSE 2 0 3 0 Black White
% 4 3320 1740 5850 2150 780 1760 5670 1780
% 
{\color[named]{Black}{%
\special{pn 8}%
\special{ar 3320 1740 2530 410  0.1047219  3.0932053}%
}}%
% ELLIPSE 2 0 3 0 Black White
% 4 3320 1740 5850 2150 780 1760 5670 1780
% 
{\color[named]{Black}{%
\special{pn 8}%
\special{ar 3320 1740 2530 410  0.1047219  3.0932053}%
}}%
% ELLIPSE 2 0 3 0 Black White
% 4 3320 3130 5850 2720 5670 3090 826 3058
% 
{\color[named]{Black}{%
\special{pn 8}%
\special{ar 3320 3130 2530 410  3.3177741  6.1784634}%
}}%
% ELLIPSE 2 0 3 0 Black White
% 4 5870 2420 5620 1770 5890 1780 5920 1800
% 
{\color[named]{Black}{%
\special{pn 8}%
\special{ar 5870 2420 250 650  4.9194617  4.7935114}%
}}%
% ELLIPSE 2 0 3 0 Black White
% 4 760 2430 510 1780 780 1790 810 1810
% 
{\color[named]{Black}{%
\special{pn 8}%
\special{ar 760 2430 250 650  4.9194617  4.7935114}%
}}%
% ELLIPSE 2 0 3 0 Black White
% 4 1890 2440 1710 2800 1890 2810 1880 2070
% 
{\color[named]{Black}{%
\special{pn 8}%
\special{ar 1890 2440 180 360  4.6583875  1.5707963}%
}}%
% ELLIPSE 2 2 3 0 Black White
% 4 1900 2440 1740 2810 1910 2010 1890 2810
% 
{\color[named]{Black}{%
\special{pn 8}%
\special{pa 1890 2810}%
\special{pa 1884 2808}%
\special{fp}%
\special{pa 1862 2800}%
\special{pa 1856 2796}%
\special{fp}%
\special{pa 1838 2782}%
\special{pa 1834 2776}%
\special{fp}%
\special{pa 1818 2758}%
\special{pa 1814 2752}%
\special{fp}%
\special{pa 1802 2732}%
\special{pa 1800 2726}%
\special{fp}%
\special{pa 1788 2704}%
\special{pa 1786 2698}%
\special{fp}%
\special{pa 1776 2674}%
\special{pa 1774 2668}%
\special{fp}%
\special{pa 1766 2644}%
\special{pa 1766 2638}%
\special{fp}%
\special{pa 1760 2614}%
\special{pa 1758 2606}%
\special{fp}%
\special{pa 1752 2582}%
\special{pa 1752 2576}%
\special{fp}%
\special{pa 1748 2550}%
\special{pa 1746 2542}%
\special{fp}%
\special{pa 1744 2518}%
\special{pa 1744 2510}%
\special{fp}%
\special{pa 1742 2484}%
\special{pa 1742 2476}%
\special{fp}%
\special{pa 1740 2450}%
\special{pa 1740 2444}%
\special{fp}%
\special{pa 1740 2416}%
\special{pa 1742 2410}%
\special{fp}%
\special{pa 1742 2384}%
\special{pa 1742 2376}%
\special{fp}%
\special{pa 1746 2350}%
\special{pa 1746 2344}%
\special{fp}%
\special{pa 1750 2318}%
\special{pa 1750 2310}%
\special{fp}%
\special{pa 1756 2286}%
\special{pa 1756 2278}%
\special{fp}%
\special{pa 1762 2254}%
\special{pa 1764 2248}%
\special{fp}%
\special{pa 1770 2224}%
\special{pa 1772 2218}%
\special{fp}%
\special{pa 1782 2194}%
\special{pa 1784 2188}%
\special{fp}%
\special{pa 1794 2164}%
\special{pa 1796 2158}%
\special{fp}%
\special{pa 1808 2138}%
\special{pa 1812 2132}%
\special{fp}%
\special{pa 1826 2112}%
\special{pa 1830 2108}%
\special{fp}%
\special{pa 1846 2092}%
\special{pa 1852 2088}%
\special{fp}%
\special{pa 1872 2076}%
\special{pa 1878 2074}%
\special{fp}%
\special{pa 1902 2070}%
\special{pa 1910 2072}%
\special{fp}%
}}%
% VECTOR 2 0 3 0 Black White
% 2 2020 1800 1920 2000
% 
{\color[named]{Black}{%
\special{pn 8}%
\special{pa 2020 1800}%
\special{pa 1920 2000}%
\special{fp}%
\special{sh 1}%
\special{pa 1920 2000}%
\special{pa 1968 1950}%
\special{pa 1944 1952}%
\special{pa 1932 1932}%
\special{pa 1920 2000}%
\special{fp}%
}}%
% VECTOR 2 0 3 0 Black White
% 2 730 3360 5850 3450
% 
{\color[named]{Black}{%
\special{pn 8}%
\special{pa 730 3360}%
\special{pa 5850 3450}%
\special{fp}%
\special{sh 1}%
\special{pa 5850 3450}%
\special{pa 5784 3430}%
\special{pa 5798 3450}%
\special{pa 5784 3470}%
\special{pa 5850 3450}%
\special{fp}%
}}%
% STR 2 0 3 0 Black White
% 4 3220 3480 3220 3580 2 0 0 0
% $\Gamma$
\put(32.2000,-35.8000){\makebox(0,0)[lb]{$\gamma$}}%
% STR 2 0 3 0 Black White
% 4 2000 1530 2000 1630 2 0 0 0
% $\varepsilon$
\put(20.0000,-17.5000){\makebox(0,0)[lb]{$\varepsilon_t$}}%
% STR 2 0 3 0 Black White
% 4 1865 3425 1865 3475 2 0 0 0
% $\gamma_t$
\put(18.4000,-35.6000){\makebox(0,0)[lb]{$\gamma_t$}}%
% LINE 2 2 3 0 Black White
% 2 1880 2855 1885 3300
% 
{\color[named]{Black}{%
\special{pn 8}%
\special{pa 1880 2856}%
\special{pa 1886 3300}%
\special{dt 0.045}%
}}%
% DOT 0 0 3 0 Black White
% 2 1880 3380 1880 3420
% 
{\color[named]{Black}{%
\special{pn 4}%
\special{sh 1}%
\special{ar 1880 3380 16 16 0  6.28318530717959E+0000}%
\special{sh 1}%
%\special{ar 1880 3420 16 16 0  6.28318530717959E+0000}%
}}%
% VECTOR 2 0 3 0 Black White
% 2 2065 2465 2065 2445
% 
{\color[named]{Black}{%
\special{pn 8}%
\special{pa 2066 2466}%
\special{pa 2066 2446}%
\special{fp}%
\special{sh 1}%
\special{pa 2066 2446}%
\special{pa 2046 2512}%
\special{pa 2066 2498}%
\special{pa 2086 2512}%
\special{pa 2066 2446}%
\special{fp}%
}}%
\end{picture}%

\vspace{0.5cm}

Then \[
\Gamma(\varepsilon,\gamma)\in H_2(\ol{U},f^{-1}(\gamma_0)+f^{-1}(\gamma_1);\Z),\quad
\mbox{with }\quad
\partial(\Gamma(\varepsilon,\gamma))=\varepsilon_1-\varepsilon_0
%\in H_1(f^{-1}(\gamma_0)+f^{-1}(\gamma_1),\Z)
\]
where $\partial:H_2(\ol{U},f^{-1}(\gamma_0)+f^{-1}(\gamma_1);\Z)
\to H_1(f^{-1}(\gamma_0)+f^{-1}(\gamma_1),\Z)$ is the boundary map.
The homology cycle $\Gamma(\varepsilon,\gamma)$ is called a {\it Lefschetz thimble}.
Define
$
\bE(\ol{U},D;\Z)\subset H_2(\ol{U},D;\Z)
$
the subgroup generated by the Lefschetz thimbles $\Gamma(\varepsilon,\gamma)$ such that
the initial and terminal points of $\gamma$ lie in $T$ 
(hence $\partial\Gamma(\varepsilon,\gamma)\subset D$).
Define $\bE(\ol{U},\Z)$ by an exact sequence
\[
0\lra \bE(\ol{U},\Z)\lra \bE(\ol{U},D;\Z)\os{\partial}{\lra} H_1(D,\Z).
\]
Write $\bE(\ol{U},D;\Q):=\bE(\ol{U},D;\Z)\ot\Q$ etc.
\begin{prop}\label{Ev}
Assume that $f$ contains a totally degenerate semistable fiber.
Then we have
\begin{equation}\label{exp-5}
\xymatrix{
0\ar[r]&\bE(\ol{U},\Q)\ar[r]\ar[d]^\cong&\bE(\ol{U},D;\Q)\ar[r]^\partial
\ar[d]^\cong&
H_1(D,\Q)\ar@{=}[d]\ar[r]&0\\
0\ar[r]&H_2(\ol{U},\Q)_0\ar[r]&H_2(\ol{U},D;\Q)_0\ar[r]^\partial&H_1(D,\Q)\ar[r]&0
}
\end{equation}
where
$H_2(\ol{U},D;\Q)_0:= H_2(\ol{U},D;\Q)/\Image\NF^B(\ol{U})
= H_2(\ol{U},D;\Q)/H_2(e(\ol{S}),\Q)$.
\end{prop}
\begin{lem}\label{exp-lem0}
The fixed part $H^1(f^{-1}(s),\Q)^{\pi_1(S,s)}$ is trivial.
Moreover we have $H^1(\ol{U},\Q)=H^1(\ol{S},\Q)$.
In particular $H^1(X,\Q)=H^1(C,\Q)$.
\end{lem}
\begin{pf}
Let $f^{-1}(P)$ be a totally degenerate semistable fiber, and
$N$ the log monodromy on $H^1(f^{-1}(s),\Q)$ around $P$.
Then the inclusion
\[
H^1(f^{-1}(s),\Q)^{\pi_1(S,s)}=\vg(S,R^1f_*\Q)\hra \ker(N)\cong H^1(f^{-1}(P),\Q)
\]
preserves the mixed Hodge structure. The LHS is of weight one, while the RHS is of weight zero
as $f^{-1}(P)$ is totally degenerate.
Therefore the inclusion must be zero, which means
$H^1(f^{-1}(s),\Q)^{\pi_1(S,s)}=0$.
Now it is easy to show $H^1(U,\Q)=H^1(S,\Q)$ by using
the Leray spectral sequence for $f:U\to S$.
The equality 
$H^1(\ol{U},\Q)=H^1(\ol{S},\Q)$ follows from this and a commutative diagram
\[
\xymatrix{
0=H^1_D(\ol{U},\Q)\ar[r]& H^1(\ol{U},\Q)\ar[r]& H^1(U,\Q)\ar[r]& H^2_D(\ol{U},\Q)
\cong H_2(D,\Q)\\
0=H^1_{T}(\ol{S},\Q)\ar[r]\ar[u]& H^1(\ol{S},\Q)\ar[r]\ar[u]& H^1(S,\Q)
\ar[r]\ar[u]^\cong& H^2_{T}(\ol{S},\Q)\cong H_0(T,\Q)\ar[u]^\cup
}\]
\end{pf}

\begin{lem}\label{Ev-lem}
The sequence
\begin{equation}\label{Ev-lem-2}
H_2(U,\Q)\lra H_2(\ol{U},D;\Q)_0\os{\partial}{\lra} H_1(D,\Q)\lra 0
\end{equation}
is exact. 
\end{lem}
\begin{pf}
The surjectivity of $\partial$ is immediate from the fact that
the composition $H_1(D,\Q)\to H_1(\ol{U},\Q)\os{\cong}{\to} H_1(\ol{S},\Q)$ is zero.
Let us show
\begin{equation}\label{Ev-lem-1}
\Image(H_2(U,\Q)\lra H_2(\ol{U},D;\Q)_0)
=\Image(H_2(\ol{U},\Q)\lra H_2(\ol{U},D;\Q)_0).
\end{equation}
Write 
\[
H^2(D,\Q)_0:=\Coker[H_2(e(\ol{S}),\Q)\to H_2^D(\ol{U},\Q)\cong H^2(D,\Q)]
\]
Consider a diagram 
\[
\xymatrix{
&H_2(D,\Q)\ar[d]^a\\
H_2(U,\Q)\ar[r]& H_2(\ol{U},\Q)_0\ar[d]\ar[r]^b& H^2(D,\Q)_0\ar[r]
&H_1(U,\Q)\ar[r]^c&H_1(\ol{U},\Q)\\
&H_2(\ol{U},D;\Q)_0\\
}\]
with exact row and column.
Hence it is enough to show $\Image(ba)=\Image(b)$ or equivalently $\dim\Coker(ba)=\dim\Coker(b)
(=\dim\ker(c))$.
Since $ba$ is given by the intersection pairing, Zariski's lemma (\cite{barth} III (8.2))
shows that
$\dim\Coker(ba)=\dim H_0(T)$ if $\ol{S}\ne C$ and $=\dim H_0(T)-1$ if $\ol{S}= C$.
On the other hand, \[\ker[H_1(U,\Q)\os{c}{\to} H_1(\ol{U},\Q)]\cong
\ker[H_1(S,\Q)\to H_1(\ol{S},\Q)]\cong \Coker[H_2(\ol{S})\to H^0(T)]\]
where we used Lemma \ref{exp-lem0} in the first isomorphism.
So we are done.
\end{pf}
\begin{lem}\label{SS-lem}
Let $f^{-1}(P)$ be a totally degenerate semistable fiber.
Let $\Van\subset H_1(f^{-1}(s),\Q)$ be the subspace generated by the vanishing cycles as
$s\to P$. Then we have \[\Q[\pi_1(S,s)](\Van)= H_1(f^{-1}(s),\Q).\]
\end{lem}
\begin{pf}
Put $V=\Q[\pi_1(S,s)](\Van)$. 
By Deligne's semisimplicity theorem (\cite{HodgeII} 4.2.6)
there is an complementary space $V'\subset H_1(f^{-1}(s),\Q)$ which is 
stable under the action of $\pi_1(S,s)$.
Let $N$ be the log monodromy around $P$. Since $\Image(N)=\Van$
one has $NV'\subset V'\cap \Van=0$.
On the other hand the composition of maps $V'\hra H_1(f^{-1}(s),\Q)/\Van\os{N}{\lra}\Van$
is injective and its image is $NV'$.
Therefore we have $V'=0$.
\end{pf}
\noindent{\it Proof of Prop.\ref{Ev}.}
Let $\cL$ be the local system on $S^o(\C)$ whose fiber is $H_1(f^{-1}(s),\Q)$.
Then the image of $H_2(U,\Q)$ in $H_2(\ol{U},D;\Q)_0$ coincides with that of
$H_1(S,\cL)$. The homology group $H_1(S,\cL)$ is generated by Lefschetz thimbles
$\Gamma(\varepsilon,\gamma)$ such that the initial and terminal points of 
$\gamma$ are the same in $S$ and $\partial\Gamma(\varepsilon,\gamma)=0$. 
Take an arbitrary path $\delta$ such that the initial point lies in $T$
and the terminal point is that of $\gamma$.
Put $\wt{\gamma}=\delta\cdot \gamma\cdot\delta^{-1}$.
Then $\Gamma(\varepsilon,\wt{\gamma})\in \bE(\ol{U},D;\Z)$ and 
the image of it in $H_2(\ol{U},D;\Q)$ coincides with that of $\Gamma(\varepsilon,\gamma)$.
This means that there is some subgroup $\bE(\ol{U},D;\Z)'\subset \bE(\ol{U},D;\Z)$
such that the image of $\bE(\ol{U},D;\Q)'$ in $H_2(\ol{U},D;\Q)$ coincides with that
of $H_2(U,\Q)$.

Next we show that the boundary map $\partial:\bE(\ol{U},D;\Q)\to H_1(D,\Q)$
is surjective.
Let $f^{-1}(P)$ be a totally degenerate semistable fiber and
$\Van\subset H_1(f^{-1}(s),\Q)$ the space of the vanishing cycles.
By Lemma \ref{SS-lem}, for any $\nu\in H_1(f^{-1}(s),\Q)$ 
there is a sum of Lefschetz thimbles $\Gamma=\sum \Gamma(\varepsilon,\gamma)$
with $\gamma\in \pi_1(S,s)$ such that 
$\partial\Gamma=\nu-$(vanishing cycle). 
By adding a path from $s$ to a point $s_0\in T$ and
a path from $s$ to a point $P$ to $\Gamma$, 
one has
a thimble $\Gamma'\in \bE(\ol{U},D;\Q)$ such that
$\partial\Gamma'=\nu\in H_1(f^{-1}(s_0),\Q)$.
This means that $\bE(\ol{U},D;\Q)\to H_1(D,\Q)$ is surjective.

There remains to show the injectivity of $\bE(\ol{U},D;\Q)\to H_2(\ol{U},D;\Q)_0$.
This is trivial unless $\ol{S}=C$. In case $\ol{S}=C$, this follows from the following
fact.
The composition
$\bE(\ol{U},D;\Q)\to H_2(\ol{U},D;\Q){\to} H_2(C,T;\Q)\cong\Q$ is zero,
while the composition 
$H_2(e(C))\to H_2(\ol{U},D;\Q){\to} H_2(C,T;\Q)\cong \Q$ is bijective.
Q.E.D.

%%%%%%%%%%%%%%%%%%%%%%%%%

\section{A formula for Regulator on $K_1$ of a fibration of curves}\label{exp-sect}
The following is the main theorem of this paper.
\begin{thm}\label{ExtMHS}
Let $f:X\to C$ be a fibration of curves over $\C$ as in \S \ref{deRh-sect-1}.
Suppose that $f$ has a totally degenerate semistable fiber.
Let $C^o\subset C$ be the maximal open set such that $\ol{\nabla}$ \eqref{PF-rem-eq}
is bijective on $C^o$. Put $X^o=f^{-1}(C^o)$.
Let 
\[
\Phi:\bE(X^o,\Q)\cong H_2(X^o,\Q)/\NF^B(X^o)\lra \Hom(\goodxa,\C)
\]
be the map of period integral defined by
\[
\Phi(\Delta)=\left[\omega\longmapsto\int_\Delta \omega\right],\quad \omega
\in \goodxa\footnote{Note that, since $\goodxa\subset
\vg(X^o,\Omega^2_{X^o})$ \eqref{hodge-pf6}, 
one can {\it a priori} define ``$\int_\Delta\omega$" 
only for $\Delta\in H_2(X^o,\Q)$.}
\]
where $\goodxa$ and $\bE(X^o,\Q)$ are as in \S \ref{good-sect} and
\S \ref{exp-sect2} respectively, and the isomorphism is due to Prop.\ref{Ev}.
\begin{enumerate}
\item[{\rm (1)}]
There is an isomorphism
\[
\Ext^1(\Q,H^2(X,\Q(2))/\NF(X))\cong\Coker(\Phi).
\]
\item[{\rm (2)}]
Let $D=\sum_i f^{-1}(P_i)$ be a union of singular fibers which are contained in $X^o$. 
Let $\xi\in H^3_{\cM,D}(X,\Q(2))$ be an arbitrary element, and put
$\gamma:=\reg_D(\xi)\in H_1^B(D,\Q)$.
Fix $\Gamma\in \bE(X^o,D;\Q)$
such that $\partial(\Gamma)=\gamma$.
Then 
\[
\reg(\xi)=
\left[\omega\longmapsto
\int_\Gamma \omega\right]\in \Coker(\Phi).
\]
\end{enumerate}
\end{thm}
\begin{pf}
Since $\goodxa\cong F^1H^2_\dR(X)_0$ by Prop. \ref{good-prop2}, 
the period map $\Phi$ factors through $H_2(X,\Q)/\NF^B(X)$.
Then (1) follows from 
the fact that $H_2(X^o,\Q)\to H_2(X,\Q)/\NF^B(X)$ is surjective.
We show (2).
Let $\delta$ be the composition of maps
\begin{align}
\Hom_\MHS(\Q,H_1(D,\Q))&\lra
\Ext^1_\MHS(\Q,H_2(X,\Q)/H_2(D,\Q))\label{exp-ext2}\\
&\lra
\Ext^1_\MHS(\Q,H_2(X,\Q)/\NF^B(X))\\
&\cong\quad\Ext^1_\MHS(\Q,H^2(X,\Q(2))/\NF(X))\label{exp-ext1}
\end{align}
where \eqref{exp-ext2} is the connecting homomorphism arising from the exact
sequence
\[
0\lra H_2(X,\Q)/H_2(D,\Q)\lra H_2(X,D;\Q)\lra H_1(D,\Q)\lra 0.
\]
Then one has $\reg(\xi)=\delta(\reg_D(\xi))=\delta(\gamma)$ by Theorem \ref{bregin-thm}.
To compute $\delta(\gamma)$, we consider a commutative diagram
\begin{equation}\label{main-diagram}
\xymatrix{
0\ar[r]&H_2(X,\Q)/H_2(D,\Q)\ar[d]_{\mbox{surj.}}\ar[r]\ar@{}[rd]|{\square}
&H_2(X,D;\Q)\ar[d]_{\mbox{surj.}}\ar[r]&H_1(D,\Q)\ar@{=}[d]\ar[r]&0\\
0\ar[r]&H_2(X,\Q)/\NF^B(X)\ar[r]&\ol{H_2(X,D;\Q)}\ar[r]&H_1(D,\Q)\ar[r]&0\\
0\ar[r]&H_2(X^o,\Q)/H_2(D)\ar[u]^{\mbox{surj.}}\ar[r]
&H_2(X^o,D;\Q)\ar[u]^{\mbox{surj.}}\ar[r]^{\quad\partial}&
H_1(D,\Q)\ar[r]\ar@{=}[u]&0.
}
\end{equation}
Then it is enough to describe the extension data of the bottom row in \eqref{main-diagram}.
For $\omega\in \goodxa$, let $\omega_{X^o,D}\in F^1H^2_\dR(X^o,D)$ denotes a lifting.
Then we have
\[
\reg(\xi)=
\left[\omega\longmapsto\int_\Gamma \omega_{X^o,D}\right]
\in \Hom(\goodxa,\C)/\Image \Phi
\]
by \eqref{beireg-3}.
On the other hand, by Theorem \ref{hodge-mot} and \eqref{pairingVD2}, we have
\[
\int_\Gamma \omega_{X^o,D}=\int_\Gamma \omega.
\]
This competes the proof.
\end{pf}
%%%%%%%%%%%%%%%%%%%%%%%%%%%%%%%%%%%%%%%%%%%%%%%
\section{Example : Elliptic fibration with $\mu_l$-action}\label{Example-sect}
%\subsection{Notation}
Let $F\subset\R$ be a subfield. We consider two polynomials 
$g_2(t),g_3(t)\in F[t]$ which satisfy the 
following (however see Remark \ref{unf-rem}).
\begin{enumerate}
\item[\bf(E1)]
$\Delta:=g_2^3-27g_3^2=ct^a(1-t)^b$ \quad for some $a,b\in\Z_{\geq 1}$ and $c\in \R_{>0}$,
\item[\bf(E2)]
$2g_2g'_3-3g'_2g_3=c't^{a'}(1-t)^{b'}$
\quad for some $a',b'\in\Z_{\geq 0}$ and $c'\ne0$,
\item[\bf(E3)]
$g_2(0),~g_2(1)>0$ and $g_3(0)g_3(1)<0$.
\item[\bf(E4)]
$g_2(t)\geq 0$ for $0\leq t\leq 1$.
\end{enumerate}
Let $l\geq 1$ and $\kappa\in F^\times$. We discuss an elliptic fibration
\begin{equation}\label{ge-0}
f:X=X_l\lra \P^1,\quad f^{-1}(t):\kappa y^2=4x^3-g_2(t^l)x-g_3(t^l)
\end{equation}
defined over $F$.
In what follows, we take $X$ to be minimal, i.e. there is no exceptional curve in a fiber. 
There is the section $e:\P^1\to X$ of ``infinity".
Let $\zeta_l$ be a $l$-th root of unity, and 
$\sigma$ an automorphism of $X_\C=X\times_F\C$ given by $(x,y,t)\mapsto (x,y,\zeta_lt)$.
Put $D:=f^{-1}(1)$ a multiplicative fiber of type $I_b$.
Let us choose $\kappa$ such that $D=f^{-1}(1)$ is split multiplicative over $F$,
or equivalently
\[
\sqrt{-6\kappa g_3(1)}\in F^\times.
\]
Then one has an element $\xi\in H^3_{\cM,D}(X,\Q(2))$ such that
$\gamma:=\reg_D(\xi)\in H^B_1(D(\C),\Q)$ is a generator (cf. Prop. \eqref{bei-hodge}).
This is uniquely determined modulo the decomposable part.
We are going to compute the real regulator
\[
\reg_\R(\xi)\in \Ext^1_{\R\mbox{-}\MHS}(\R,H^2(X)_\ind\ot\R(1))^{F_\infty=1},\quad
H^2(X)_\ind:=H^2(X,\Q(1))/\NS(X)
\]
where $F_\infty$ is the infinite Frobenius.

\subsection{Computation of $\Lambda(X)_\mathrm{rat}$}\label{good-expm-sect}
The elliptic fibration \eqref{ge-0} 
is smooth outside $t=0,\infty,\zeta_l^i$, $(0\leq i \leq l-1)$ by {\bf (E1)}.
It follows from {\bf(E3)} that
$f^{-1}(0)$ is a semistable fiber of type $I_{al}$, and
$f^{-1}(\zeta_l^i)$ is of type $I_b$ (cf. Tate's algorithm, \cite{Si} IV).
Let $\nu_\infty$ be the number of irreducible components of $f^{-1}(\infty)$,
and $\varepsilon_\infty$ the Kodaira index of $f^{-1}(\infty)$:
\begin{center}
\begin{tabular}{c|ccccccccc}
$\varepsilon_\infty$&0&$b$&2&3&4&$b+6$&10&9&8\\
\hline
$f^{-1}(\infty)$&smooth&$\mathrm{I}_b$&II&III&IV&$\mathrm{I}^*_b$&II*&III*&IV*
\end{tabular}
\end{center}
Note that
\begin{equation}\label{ge-1}
\varepsilon_\infty=
\begin{cases}
\nu_\infty-1&\mbox{smooth}\\
\nu_\infty&\mbox{multiplicative}\\
\nu_\infty+1&\mbox{additive}.
\end{cases}
\end{equation}
As is well-known, we have
\begin{equation}\label{ge-2}
\begin{cases}
K_X\cong f^*\O_{\P^1}(\frac{al+bl+\varepsilon_\infty}{12}-2) \quad
\mbox{(Kodaira's canonical bundle formula)}\\
h^{20}(X):=\dim H^0(\Omega^2_X)=\frac{al+bl+\varepsilon_\infty}{12}-1\\
b^2:=\mathrm{rank} H^2(X_\C)=al+bl+\varepsilon_\infty-2\\
\rho_f:=\mathrm{rank} \NF(X_\C)=al+(b-1)l+\nu_\infty.
\end{cases}
\end{equation}

\begin{lem}\label{ge-3}
Let $s=t^{-1}$ and $k\geq 0$ be the minimal integer such that
both of $\bar{g}_2(s):=s^{4k}g_2(s^{-l})$ and $\bar{g}_3(s):=s^{6k}g_3(s^{-l})$ have no pole.
This is equivalent to saying that $k$ is the integer such that
\[
ay_1^2=4x_1^3-\bar{g}_2(s)x_1-\bar{g}_3(s),\quad x_1:=s^{2k}x,~y_1=s^{3k}y
\]
is the minimal Weierstrass equation of $X$ over $s=0$ $(t=\infty)$.
Then \[
\vg(X,\Omega_X^2)=\langle t^{i-1}dt\frac{dx}{y}~|~1\leq i \leq k-1\rangle.
\]
In particular, $h^{20}=k-1$.
\end{lem}
\begin{pf}
$f_*K_X\cong\O_{\P^1}(\frac{al+bl+\varepsilon_\infty}{12}-2)$ 
is a locally free sheaf of rank one.
This has a free basis $dtdx/y$ on $\P^1\setminus\{\infty\}$ and
$dsdx_1/y_1$ on a neighborhood of $s=0$ ($t=\infty$).
Then the assertion follows from
\[
t^{i-1}dt\frac{dx}{y}=-s^{k-i-1}ds\frac{dx_1}{y_1}.
\]
\end{pf}
\begin{prop}\label{NS=NF}
Suppose that $l$ is a prime number and $h^{20}(X)>0$.
Then $\dim H^2(X)_\ind=l-1$ and $\NF(X_\C)\ot\Q=\NS(X_\C)\ot\Q$
(hence $H^2(X)_0\cong H^2(X)_\ind$).
Moreover $f^{-1}(\infty)$ is an additive fiber.
\end{prop}
\begin{pf}
By \eqref{ge-1} and \eqref{ge-2}, we have
$b^2-\rho_f=l-2+(\varepsilon_\infty-\nu_\infty)\leq l-1$.
On the other hand, $\sigma$ acts on $H^2(X_\C,\Q)/\NF(X_\C)$ and it has an eigenvalue
$\zeta_l$ since $dtdx/y\in \vg(X,\Omega^2_X)$ by Lemma \ref{ge-3}. Since $l$ is a prime number,
the characteristic polynomial of $\sigma$ must be divided by $1+x+x^2+\cdots+x^{l-1}$,
and hence its degree is at least $l-1$. This implies $b^2-\rho_f\geq l-1$.
Hence we have $b^2-\rho_f=l-1$ and $\varepsilon_\infty-\nu_\infty=1$. 
This implies that $f^{-1}(\infty)$ is an additive fiber by \eqref{ge-1}.
Let $\rho:=\mathrm{rank}\NS(X_\C)$. Obviously $\rho\geq \rho_f$.
Since $\sigma$ acts on $H^2(X_\C,\Q)/\NS(X_\C)$ as well, the same argument yields
$b^2-\rho\geq l-1$.
We thus have $\rho\leq \rho_f$ and hence $\rho=\rho_f$.
\end{pf}
\begin{prop}\label{goodbasis}
Suppose that $l$ is a prime number and $h^{20}(X)>0$.
Then
\[
\goodxb=\langle t^{i-1}dt\frac{dx}{y}~|~1\leq i \leq h^{20}\rangle\cong F^2H_\dR(X),
\]
\[
\goodxa=\langle t^{i-1}dt\frac{dx}{y}~|~1\leq i \leq l-1-h^{20}\rangle\cong
F^1H_\dR(X)_\ind.
\]
\end{prop}
\begin{pf}
The former was shown in Lemma \ref{ge-3}.
We show the latter. Let $s=t^{-1}$ and $k$, $x_1=s^{2k}x$ and $y_1=s^{3k}y$
be as in Lemma \ref{ge-3}.
Put $T:=\{0,\infty,\zeta_l^i~|~0\leq i \leq l-1\}\subset \P^1$ and
$U:=f^{-1}(S)\to S:=\P^1-T$.
Let $\cH_e$ be Deligne's canonical extension of $\cH=R^1f_*\Omega^\bullet_{U/S}$.
Let $\omega,\omega^*\in \vg(\P^1\setminus T,\cH)$ 
be as in Lemma \ref{he-basis}. They give a free basis of $\cH_e$ 
over $\P^1\setminus\{\infty\}$ by Theorem \ref{A2}.
Since $f^{-1}(\infty)$ is an additive fiber, 
$\{t^{k-1}\omega,t^{-k}\omega^*\}$ is a free basis on a neighborhood
of $\infty$. We thus have
\begin{equation}\label{gb-1}
\cH_e^{1,0}\cong\O_{\P^1}(k-1),\quad
\cH_e^{0,1}\cong\O_{\P^1}(-k),
\end{equation}
\begin{equation}\label{gb-2}
\vg(\P^1,\Omega^1_{\P^1}(\log T)\ot\cH_e)=\left\langle\frac{t^{i-1}dt}{t(1-t^l)}\ot\omega,~
\frac{t^{j-1}dt}{t(1-t^l)}\ot\omega^*~|~
\begin{matrix}1\leq i\leq l+k,\\
1\leq j\leq l-k+1
\end{matrix}
\right\rangle
\end{equation}
\begin{equation}\label{gb-3}
\vg(\P^1,\cH_e^{1,0})=\langle t^{i-1}\frac{dx}{y}~|~1\leq i \leq k\rangle
\end{equation}
and 
\[
H^1_\zar(\P^1,\cH^{1,0}_e\to \Omega^1_{\P^1}(\log T)\ot\cH_e)\cong\Coker[
\vg(\P^1,\cH^{1,0}_e)\to \vg(\P^1,\Omega^1_{\P^1}(\log T)\ot\cH_e)].
\]
By Theorem \ref{A1} \eqref{A1-3} we have
\begin{equation}\label{gb-4}
\cH_e^{1,0}\lra \Omega^1_{\P^1}(\log T)\ot\cH_e^{0,1},
\quad \omega\longmapsto \frac{3(2g_2g'_3-3g'_2g_3)}{4\Delta}dt\ot\omega^*.
\end{equation}
By {\bf (E1)} and {\bf (E2)} we have  
\[
\frac{1}{27}\frac{g_2}{g_3}\frac{dj}{j}=\frac{2g_2g'_3-3g'_2g_3}{\Delta}dt=
\frac{dt}{t^{a-a'}(1-t)^{b-b'}}\times(\mbox{const.})
,\quad j:=\frac{1728g_2^3}{\Delta}.
\]
By {\bf (E3)}, we have $a-a'=b-b'=1$:
\begin{equation}\label{gb-5}
\frac{3(2g_2g'_3-3g'_2g_3)}{4\Delta}dt=\frac{dt}{t(1-t^l)}\times\mbox{(const.)}.
\end{equation}
This shows that \eqref{gb-4} is bijective on $\P^1\setminus\{\infty\}$.
Let 
\[
\phi:H^1_\zar(\P^1,\cH^{1,0}_e\to \Omega^1_{\P^1}(\log T)\ot\cH_e)
\lra \vg(S,\Omega^1_S\ot\cH)
\]
be the composition of \eqref{hodge-pf2-0} and \eqref{hodge-pf2-1}.
Then by \eqref{gb-5} and \eqref{A1-3}, we have
\begin{align*}
\phi\left(\frac{t^{j-1}dt}{t(1-t^l)}\ot\omega^*\right)
&=\left(-t^{j-1}\frac{d\Delta}{12\Delta}+(j-1)t^{j-2}dt\right)\frac{dx}{y}\times\mbox{(const.)}\\
&=h(t)\times\frac{dt}{t(1-t^l)}\frac{dx}{y}
\end{align*}
with deg $h(t)\leq l+j-1$. Hence
\[
\goodua:=\Image(\phi)\subset \left\langle\frac{t^{j-1}dt}{t(1-t^l)}\frac{dx}{y}~|~
1\leq j\leq 2l-k+1\right\rangle.
\]
This yields
\[
\goodxa\subset \left\langle t^{j-1}dt\frac{dx}{y}~|~
1\leq j\leq l-k=l-1-h^{20}\right\rangle.
\]
Since $\dim \goodxa=\dim F^1H^2_\dR(X)_0=l-1-h^{20}$ 
by Prop. \ref{NS=NF}, the equality holds in the above. This is the desired assertion.
\end{pf}

\begin{rem}\label{unf-rem}
Since $\dim H^2(X)_\ind\cap H^{1,1}=l-1-2h^{2,0}\geq0$, one has
$l-1-((al+bl+\varepsilon_\infty)/6-2)\geq 0$ for any large prime number $l$.
This implies $a+b\leq 6$, together with $a'=a-1$ and $b'=b-1$ by \eqref{gb-5}.
Then by case-by-case analysis based on $(a,b)$, one can show that there are only
the following pairs of $(g_2,g_3)$ satisfying {\bf(E1)},\ldots,{\bf(E4)}, up to the
equivalence $(g_2,g_3)\sim(h^4g_2,h^6g_3)$ or $(g_2(t),g_3(t))\sim(g_2(1-t),g_3(1-t))$.
\begin{enumerate}
\item[(i)]
$(g_2,g_3)=(3,1-2t)$, $(a,b)=(1,1)$
\item[(ii)]
$(g_2,g_3)=(12-9t,8-9t)$, $(a,b)=(2,1)$
\item[(iii)]
$(g_2,g_3)=(27-24t,-8t^2+36t-27)$, $(a,b)=(3,1)$
\item[(iv)]
$(g_2,g_3)=(3(t^2-16t+16),(t-2)(t^2+32t-32))$, $(a,b)=(4,1)$
\item[(v)]
$(g_2,g_3)=(12(t^2-t+1),4(t-2)(t+1)(2t-1))$, $(a,b)=(2,2)$.
\end{enumerate}
\end{rem}
\subsection{Computation of Lefschetz thimbles : Cycles $\Delta$ and $\Gamma$}
Let $\delta_0$ (resp. $\delta_1$) be the homology cycle in $H_1(f^{-1}(t),\Z)$
which vanishes as $t\to 0$ (resp. $t\to 1$).
Define $\Delta$ and $\Gamma$ to be fibrations
over the segment $[0,1]\subset \P^1(\C)$ whose fibers are the vanishing cycles
$\delta_1$ and $\delta_0$ respectively.
\[
\Delta\in H_2(\ol{U},f^{-1}(0);\Z),\quad
\Gamma\in H_2(\ol{U},f^{-1}(1);\Z).
\]
The boundary $\partial\Delta$ (resp. $\partial\Gamma$) is a generator
of the homology group $H_1(f^{-1}(0),\Z)$ (resp. $H_1(f^{-1}(1),\Z)$).

%%%%%%%%%%%%%%%%%%%%%%%%%%%%%
\vspace{0.5cm}

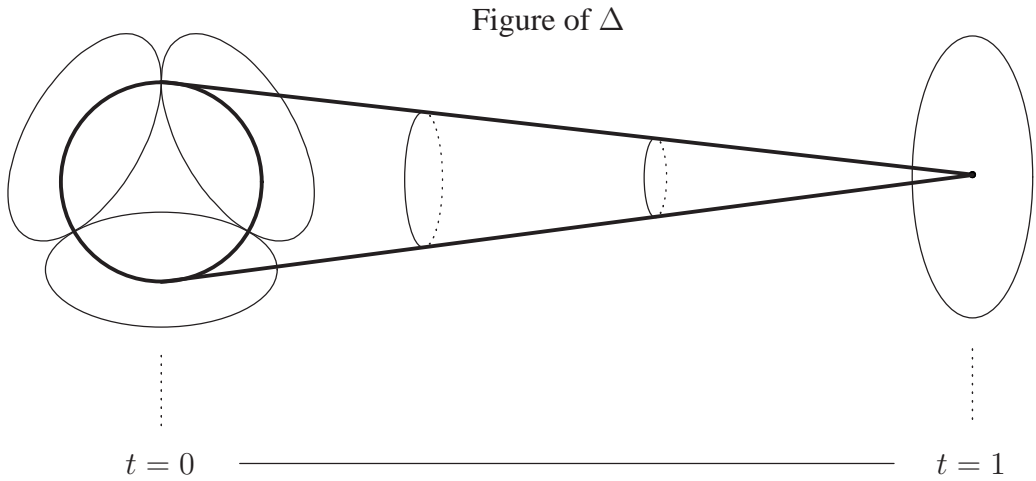
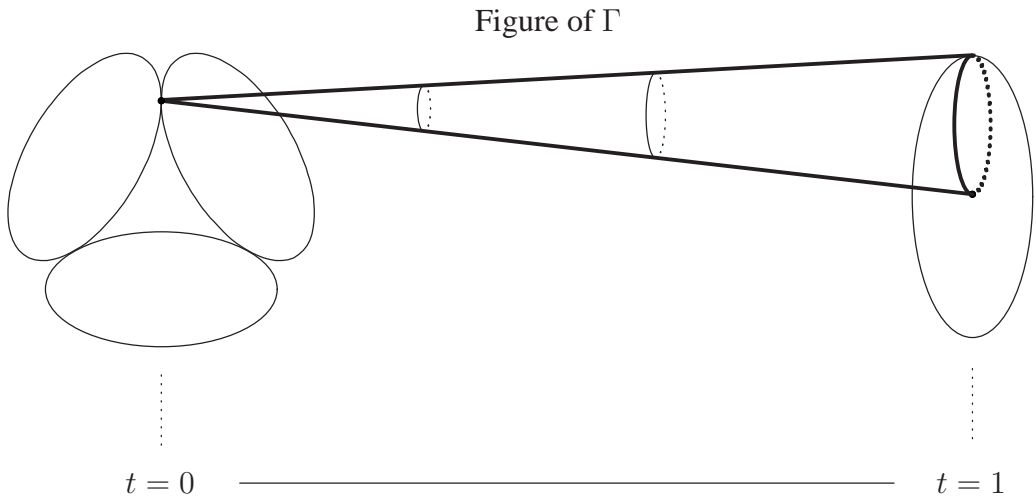
\begin{center}
%WinTpicVersion4.23
\unitlength 0.1in
\begin{picture}( 53.0500, 53.6500)(  9.0200,-62.0000)
% LINE 2 2 3 0 Black White
% 4 5896 2602 5896 3002 1696 3002 1696 2602
% 
{\color[named]{Black}{%
\special{pn 8}%
\special{pa 5896 2602}%
\special{pa 5896 3002}%
\special{dt 0.045}%
\special{pa 1696 3002}%
\special{pa 1696 2602}%
\special{dt 0.045}%
}}%
% ELLIPSE 2 0 3 0 Black White
% 4 5896 4704 6207 5440 6207 5440 6207 5440
% 
{\color[named]{Black}{%
\special{pn 8}%
\special{ar 5896 4704 312 736  0.0000000  6.2831853}%
}}%
% DOT 2 0 3 0 Black White
% 2 5896 4692 5896 4692
% 
{\color[named]{Black}{%
\special{pn 4}%
\special{sh 1}%
\special{ar 5896 4692 8 8 0  6.28318530717959E+0000}%
\special{sh 1}%
\special{ar 5896 4692 8 8 0  6.28318530717959E+0000}%
}}%
% LINE 2 2 3 0 Black White
% 4 5896 5602 5896 6002 1696 6002 1696 5602
% 
{\color[named]{Black}{%
\special{pn 8}%
\special{pa 5896 5602}%
\special{pa 5896 6002}%
\special{dt 0.045}%
\special{pa 1696 6002}%
\special{pa 1696 5602}%
\special{dt 0.045}%
}}%
% STR 2 0 3 0 Black White
% 4 1690 6100 1690 6200 5 0 0 0
% $t=0$
\put(16.9000,-62.0000){\makebox(0,0){$t=0$}}%
% STR 2 0 3 0 Black White
% 4 5890 6100 5890 6200 5 0 0 0
% $t=1$
\put(58.9000,-62.0000){\makebox(0,0){$t=1$}}%
% LINE 2 0 3 0 Black White
% 2 5490 6200 2102 6200
% 
{\color[named]{Black}{%
\special{pn 8}%
\special{pa 5490 6200}%
\special{pa 2102 6200}%
\special{fp}%
}}%
% ELLIPSE 2 0 3 0 Black White
% 4 1696 5187 2296 4887 2296 4887 2296 4887
% 
{\color[named]{Black}{%
\special{pn 8}%
\special{ar 1696 5188 600 300  0.0000000  6.2831853}%
}}%
% SPLINE 2 0 3 0 Black White
% 41 1598 3978 1635 4007 1664 4049 1684 4102 1693 4165 1694 4236 1685 4313 1665 4395 1638 4480 1602 4564 1558 4647 1508 4727 1453 4800 1393 4866 1332 4925 1270 4971 1208 5005 1148 5029 1093 5038 1042 5034 998 5017 961 4987 932 4945 913 4893 903 4829 903 4759 911 4682 931 4599 959 4515 995 4431 1038 4347 1089 4268 1144 4195 1203 4128 1264 4070 1327 4024 1388 3989 1448 3966 1503 3956 1554 3960 1598 3978
% 
{\color[named]{Black}{%
\special{pn 8}%
\special{pa 1598 3978}%
\special{pa 1624 3998}%
\special{pa 1646 4020}%
\special{pa 1664 4048}%
\special{pa 1676 4078}%
\special{pa 1686 4108}%
\special{pa 1692 4140}%
\special{pa 1696 4204}%
\special{pa 1692 4268}%
\special{pa 1688 4298}%
\special{pa 1682 4330}%
\special{pa 1674 4362}%
\special{pa 1666 4392}%
\special{pa 1658 4424}%
\special{pa 1638 4484}%
\special{pa 1626 4514}%
\special{pa 1612 4542}%
\special{pa 1598 4572}%
\special{pa 1570 4628}%
\special{pa 1554 4656}%
\special{pa 1536 4684}%
\special{pa 1520 4710}%
\special{pa 1502 4738}%
\special{pa 1482 4764}%
\special{pa 1464 4788}%
\special{pa 1442 4814}%
\special{pa 1422 4836}%
\special{pa 1378 4884}%
\special{pa 1330 4928}%
\special{pa 1306 4948}%
\special{pa 1280 4966}%
\special{pa 1224 4998}%
\special{pa 1196 5012}%
\special{pa 1166 5024}%
\special{pa 1136 5032}%
\special{pa 1104 5038}%
\special{pa 1072 5038}%
\special{pa 1040 5034}%
\special{pa 1010 5024}%
\special{pa 982 5006}%
\special{pa 958 4984}%
\special{pa 940 4958}%
\special{pa 926 4930}%
\special{pa 916 4900}%
\special{pa 908 4868}%
\special{pa 904 4836}%
\special{pa 902 4804}%
\special{pa 902 4772}%
\special{pa 904 4740}%
\special{pa 912 4676}%
\special{pa 920 4646}%
\special{pa 928 4614}%
\special{pa 936 4584}%
\special{pa 956 4524}%
\special{pa 980 4464}%
\special{pa 994 4434}%
\special{pa 1022 4378}%
\special{pa 1054 4322}%
\special{pa 1072 4294}%
\special{pa 1108 4242}%
\special{pa 1128 4218}%
\special{pa 1148 4192}%
\special{pa 1168 4168}%
\special{pa 1212 4120}%
\special{pa 1234 4098}%
\special{pa 1258 4076}%
\special{pa 1282 4056}%
\special{pa 1308 4038}%
\special{pa 1336 4020}%
\special{pa 1362 4004}%
\special{pa 1390 3988}%
\special{pa 1420 3976}%
\special{pa 1450 3966}%
\special{pa 1482 3958}%
\special{pa 1514 3956}%
\special{pa 1546 3958}%
\special{pa 1576 3968}%
\special{pa 1598 3978}%
\special{fp}%
}}%
% SPLINE 2 0 3 0 Black White
% 41 1792 3978 1755 4007 1726 4049 1706 4102 1697 4165 1696 4236 1705 4313 1725 4395 1752 4480 1788 4564 1832 4647 1882 4727 1937 4800 1997 4866 2058 4925 2120 4971 2182 5005 2242 5029 2297 5038 2348 5034 2392 5017 2429 4987 2458 4945 2477 4893 2487 4829 2487 4759 2479 4682 2459 4599 2431 4515 2395 4431 2352 4347 2301 4268 2246 4195 2187 4128 2126 4070 2063 4024 2002 3989 1942 3966 1887 3956 1836 3960 1792 3978
% 
{\color[named]{Black}{%
\special{pn 8}%
\special{pa 1792 3978}%
\special{pa 1766 3998}%
\special{pa 1744 4020}%
\special{pa 1728 4048}%
\special{pa 1714 4078}%
\special{pa 1706 4108}%
\special{pa 1700 4140}%
\special{pa 1696 4204}%
\special{pa 1696 4236}%
\special{pa 1698 4268}%
\special{pa 1702 4298}%
\special{pa 1708 4330}%
\special{pa 1716 4362}%
\special{pa 1724 4392}%
\special{pa 1734 4424}%
\special{pa 1754 4484}%
\special{pa 1766 4514}%
\special{pa 1778 4542}%
\special{pa 1792 4572}%
\special{pa 1806 4600}%
\special{pa 1854 4684}%
\special{pa 1872 4710}%
\special{pa 1890 4738}%
\special{pa 1908 4764}%
\special{pa 1928 4788}%
\special{pa 1948 4814}%
\special{pa 1970 4836}%
\special{pa 2014 4884}%
\special{pa 2036 4906}%
\special{pa 2060 4928}%
\special{pa 2086 4948}%
\special{pa 2112 4966}%
\special{pa 2138 4982}%
\special{pa 2166 4998}%
\special{pa 2196 5012}%
\special{pa 2226 5024}%
\special{pa 2256 5032}%
\special{pa 2288 5038}%
\special{pa 2320 5038}%
\special{pa 2352 5034}%
\special{pa 2382 5024}%
\special{pa 2408 5006}%
\special{pa 2432 4984}%
\special{pa 2452 4958}%
\special{pa 2466 4930}%
\special{pa 2476 4900}%
\special{pa 2482 4868}%
\special{pa 2486 4836}%
\special{pa 2488 4804}%
\special{pa 2488 4772}%
\special{pa 2484 4708}%
\special{pa 2478 4676}%
\special{pa 2472 4646}%
\special{pa 2464 4614}%
\special{pa 2454 4584}%
\special{pa 2446 4554}%
\special{pa 2398 4434}%
\special{pa 2384 4406}%
\special{pa 2368 4378}%
\special{pa 2354 4350}%
\special{pa 2338 4322}%
\special{pa 2320 4294}%
\special{pa 2302 4268}%
\special{pa 2282 4242}%
\special{pa 2264 4218}%
\special{pa 2244 4192}%
\special{pa 2224 4168}%
\special{pa 2180 4120}%
\special{pa 2158 4098}%
\special{pa 2134 4076}%
\special{pa 2108 4056}%
\special{pa 2056 4020}%
\special{pa 2000 3988}%
\special{pa 1970 3976}%
\special{pa 1940 3966}%
\special{pa 1908 3958}%
\special{pa 1876 3956}%
\special{pa 1844 3958}%
\special{pa 1814 3968}%
\special{pa 1792 3978}%
\special{fp}%
}}%
% ELLIPSE 2 0 3 0 Black White
% 4 5896 1704 6207 2440 6207 2440 6207 2440
% 
{\color[named]{Black}{%
\special{pn 8}%
\special{ar 5896 1704 312 736  0.0000000  6.2831853}%
}}%
% DOT 2 0 3 0 Black White
% 2 5896 1692 5896 1692
% 
{\color[named]{Black}{%
\special{pn 4}%
\special{sh 1}%
\special{ar 5896 1692 8 8 0  6.28318530717959E+0000}%
\special{sh 1}%
\special{ar 5896 1692 8 8 0  6.28318530717959E+0000}%
}}%
% LINE 2 2 3 0 Black White
% 4 5896 2602 5896 3002 1696 3002 1696 2602
% 
{\color[named]{Black}{%
\special{pn 8}%
\special{pa 5896 2602}%
\special{pa 5896 3002}%
\special{dt 0.045}%
\special{pa 1696 3002}%
\special{pa 1696 2602}%
\special{dt 0.045}%
}}%
% ELLIPSE 2 0 3 0 Black White
% 4 1696 2187 2296 1887 2296 1887 2296 1887
% 
{\color[named]{Black}{%
\special{pn 8}%
\special{ar 1696 2188 600 300  0.0000000  6.2831853}%
}}%
% SPLINE 2 0 3 0 Black White
% 41 1598 978 1635 1007 1664 1049 1684 1102 1693 1165 1694 1236 1685 1313 1665 1395 1638 1480 1602 1564 1558 1647 1508 1727 1453 1800 1393 1866 1332 1925 1270 1971 1208 2005 1148 2029 1093 2038 1042 2034 998 2017 961 1987 932 1945 913 1893 903 1829 903 1759 911 1682 931 1599 959 1515 995 1431 1038 1347 1089 1268 1144 1195 1203 1128 1264 1070 1327 1024 1388 989 1448 966 1503 956 1554 960 1598 978
% 
{\color[named]{Black}{%
\special{pn 8}%
\special{pa 1598 978}%
\special{pa 1624 998}%
\special{pa 1646 1020}%
\special{pa 1664 1048}%
\special{pa 1676 1078}%
\special{pa 1686 1108}%
\special{pa 1692 1140}%
\special{pa 1696 1204}%
\special{pa 1692 1268}%
\special{pa 1688 1298}%
\special{pa 1682 1330}%
\special{pa 1674 1362}%
\special{pa 1666 1392}%
\special{pa 1658 1424}%
\special{pa 1638 1484}%
\special{pa 1626 1514}%
\special{pa 1612 1542}%
\special{pa 1598 1572}%
\special{pa 1570 1628}%
\special{pa 1554 1656}%
\special{pa 1536 1684}%
\special{pa 1520 1710}%
\special{pa 1502 1738}%
\special{pa 1482 1764}%
\special{pa 1464 1788}%
\special{pa 1442 1814}%
\special{pa 1422 1836}%
\special{pa 1378 1884}%
\special{pa 1330 1928}%
\special{pa 1306 1948}%
\special{pa 1280 1966}%
\special{pa 1224 1998}%
\special{pa 1196 2012}%
\special{pa 1166 2024}%
\special{pa 1136 2032}%
\special{pa 1104 2038}%
\special{pa 1072 2038}%
\special{pa 1040 2034}%
\special{pa 1010 2024}%
\special{pa 982 2006}%
\special{pa 958 1984}%
\special{pa 940 1958}%
\special{pa 926 1930}%
\special{pa 916 1900}%
\special{pa 908 1868}%
\special{pa 904 1836}%
\special{pa 902 1804}%
\special{pa 902 1772}%
\special{pa 904 1740}%
\special{pa 912 1676}%
\special{pa 920 1646}%
\special{pa 928 1614}%
\special{pa 936 1584}%
\special{pa 956 1524}%
\special{pa 980 1464}%
\special{pa 994 1434}%
\special{pa 1022 1378}%
\special{pa 1054 1322}%
\special{pa 1072 1294}%
\special{pa 1108 1242}%
\special{pa 1128 1218}%
\special{pa 1148 1192}%
\special{pa 1168 1168}%
\special{pa 1212 1120}%
\special{pa 1234 1098}%
\special{pa 1258 1076}%
\special{pa 1282 1056}%
\special{pa 1308 1038}%
\special{pa 1336 1020}%
\special{pa 1362 1004}%
\special{pa 1390 988}%
\special{pa 1420 976}%
\special{pa 1450 966}%
\special{pa 1482 958}%
\special{pa 1514 956}%
\special{pa 1546 958}%
\special{pa 1576 968}%
\special{pa 1598 978}%
\special{fp}%
}}%
% SPLINE 2 0 3 0 Black White
% 41 1792 978 1755 1007 1726 1049 1706 1102 1697 1165 1696 1236 1705 1313 1725 1395 1752 1480 1788 1564 1832 1647 1882 1727 1937 1800 1997 1866 2058 1925 2120 1971 2182 2005 2242 2029 2297 2038 2348 2034 2392 2017 2429 1987 2458 1945 2477 1893 2487 1829 2487 1759 2479 1682 2459 1599 2431 1515 2395 1431 2352 1347 2301 1268 2246 1195 2187 1128 2126 1070 2063 1024 2002 989 1942 966 1887 956 1836 960 1792 978
% 
{\color[named]{Black}{%
\special{pn 8}%
\special{pa 1792 978}%
\special{pa 1766 998}%
\special{pa 1744 1020}%
\special{pa 1728 1048}%
\special{pa 1714 1078}%
\special{pa 1706 1108}%
\special{pa 1700 1140}%
\special{pa 1696 1204}%
\special{pa 1696 1236}%
\special{pa 1698 1268}%
\special{pa 1702 1298}%
\special{pa 1708 1330}%
\special{pa 1716 1362}%
\special{pa 1724 1392}%
\special{pa 1734 1424}%
\special{pa 1754 1484}%
\special{pa 1766 1514}%
\special{pa 1778 1542}%
\special{pa 1792 1572}%
\special{pa 1806 1600}%
\special{pa 1854 1684}%
\special{pa 1872 1710}%
\special{pa 1890 1738}%
\special{pa 1908 1764}%
\special{pa 1928 1788}%
\special{pa 1948 1814}%
\special{pa 1970 1836}%
\special{pa 2014 1884}%
\special{pa 2036 1906}%
\special{pa 2060 1928}%
\special{pa 2086 1948}%
\special{pa 2112 1966}%
\special{pa 2138 1982}%
\special{pa 2166 1998}%
\special{pa 2196 2012}%
\special{pa 2226 2024}%
\special{pa 2256 2032}%
\special{pa 2288 2038}%
\special{pa 2320 2038}%
\special{pa 2352 2034}%
\special{pa 2382 2024}%
\special{pa 2408 2006}%
\special{pa 2432 1984}%
\special{pa 2452 1958}%
\special{pa 2466 1930}%
\special{pa 2476 1900}%
\special{pa 2482 1868}%
\special{pa 2486 1836}%
\special{pa 2488 1804}%
\special{pa 2488 1772}%
\special{pa 2484 1708}%
\special{pa 2478 1676}%
\special{pa 2472 1646}%
\special{pa 2464 1614}%
\special{pa 2454 1584}%
\special{pa 2446 1554}%
\special{pa 2398 1434}%
\special{pa 2384 1406}%
\special{pa 2368 1378}%
\special{pa 2354 1350}%
\special{pa 2338 1322}%
\special{pa 2320 1294}%
\special{pa 2302 1268}%
\special{pa 2282 1242}%
\special{pa 2264 1218}%
\special{pa 2244 1192}%
\special{pa 2224 1168}%
\special{pa 2180 1120}%
\special{pa 2158 1098}%
\special{pa 2134 1076}%
\special{pa 2108 1056}%
\special{pa 2056 1020}%
\special{pa 2000 988}%
\special{pa 1970 976}%
\special{pa 1940 966}%
\special{pa 1908 958}%
\special{pa 1876 956}%
\special{pa 1844 958}%
\special{pa 1814 968}%
\special{pa 1792 978}%
\special{fp}%
}}%
% CIRCLE 0 0 3 0 Black White
% 4 1696 1729 1696 2249 1696 2249 1696 2249
% 
{\color[named]{Black}{%
\special{pn 20}%
\special{ar 1696 1730 520 520  0.0000000  6.2831853}%
}}%
% LINE 0 0 3 0 Black White
% 2 1694 1208 5896 1692
% 
{\color[named]{Black}{%
\special{pn 20}%
\special{pa 1694 1208}%
\special{pa 5896 1692}%
\special{fp}%
}}%
% DOT 0 0 3 0 Black White
% 2 5896 1692 5896 1692
% 
{\color[named]{Black}{%
\special{pn 4}%
\special{sh 1}%
\special{ar 5896 1692 16 16 0  6.28318530717959E+0000}%
\special{sh 1}%
\special{ar 5896 1692 16 16 0  6.28318530717959E+0000}%
}}%
% LINE 0 0 3 0 Black White
% 2 1696 2252 5896 1694
% 
{\color[named]{Black}{%
\special{pn 20}%
\special{pa 1696 2252}%
\special{pa 5896 1694}%
\special{fp}%
}}%
% ELLIPSE 2 0 3 0 Black White
% 4 3054 1718 2956 1364 3054 1366 3054 2074
% 
{\color[named]{Black}{%
\special{pn 8}%
\special{ar 3054 1718 98 354  1.5707963  4.7123890}%
}}%
% ELLIPSE 2 2 3 0 Black White
% 4 3054 1718 3152 1364 3054 2074 3054 1366
% 
{\color[named]{Black}{%
\special{pn 8}%
\special{pa 3054 1364}%
\special{pa 3059 1364}%
\special{pa 3059 1365}%
\special{pa 3061 1365}%
\special{fp}%
\special{pa 3086 1383}%
\special{pa 3086 1384}%
\special{pa 3087 1384}%
\special{pa 3087 1385}%
\special{pa 3089 1387}%
\special{pa 3089 1388}%
\special{pa 3090 1388}%
\special{pa 3090 1388}%
\special{fp}%
\special{pa 3105 1416}%
\special{pa 3106 1417}%
\special{pa 3106 1419}%
\special{pa 3107 1419}%
\special{pa 3107 1421}%
\special{pa 3108 1422}%
\special{pa 3108 1423}%
\special{fp}%
\special{pa 3119 1454}%
\special{pa 3119 1454}%
\special{pa 3120 1455}%
\special{pa 3120 1457}%
\special{pa 3121 1458}%
\special{pa 3121 1461}%
\special{fp}%
\special{pa 3130 1494}%
\special{pa 3130 1496}%
\special{pa 3131 1497}%
\special{pa 3131 1501}%
\special{pa 3131 1501}%
\special{fp}%
\special{pa 3138 1536}%
\special{pa 3138 1538}%
\special{pa 3139 1539}%
\special{pa 3139 1543}%
\special{fp}%
\special{pa 3144 1579}%
\special{pa 3144 1582}%
\special{pa 3145 1583}%
\special{pa 3145 1587}%
\special{fp}%
\special{pa 3149 1624}%
\special{pa 3149 1624}%
\special{pa 3149 1631}%
\special{fp}%
\special{pa 3151 1668}%
\special{pa 3151 1676}%
\special{fp}%
\special{pa 3152 1714}%
\special{pa 3152 1722}%
\special{fp}%
\special{pa 3151 1759}%
\special{pa 3151 1767}%
\special{fp}%
\special{pa 3149 1805}%
\special{pa 3149 1811}%
\special{pa 3148 1812}%
\special{pa 3148 1812}%
\special{fp}%
\special{pa 3145 1849}%
\special{pa 3145 1853}%
\special{pa 3144 1854}%
\special{pa 3144 1857}%
\special{fp}%
\special{pa 3139 1893}%
\special{pa 3139 1896}%
\special{pa 3138 1897}%
\special{pa 3138 1900}%
\special{fp}%
\special{pa 3131 1935}%
\special{pa 3131 1938}%
\special{pa 3130 1939}%
\special{pa 3130 1943}%
\special{fp}%
\special{pa 3121 1977}%
\special{pa 3120 1978}%
\special{pa 3120 1981}%
\special{pa 3119 1981}%
\special{pa 3119 1983}%
\special{fp}%
\special{pa 3107 2015}%
\special{pa 3107 2015}%
\special{pa 3107 2016}%
\special{pa 3106 2017}%
\special{pa 3106 2019}%
\special{pa 3105 2020}%
\special{pa 3105 2021}%
\special{pa 3104 2022}%
\special{fp}%
\special{pa 3088 2050}%
\special{pa 3087 2051}%
\special{pa 3087 2052}%
\special{pa 3086 2052}%
\special{pa 3086 2053}%
\special{pa 3085 2053}%
\special{pa 3085 2054}%
\special{pa 3084 2054}%
\special{pa 3084 2054}%
\special{fp}%
\special{pa 3061 2071}%
\special{pa 3059 2071}%
\special{pa 3059 2072}%
\special{pa 3054 2072}%
\special{fp}%
}}%
% ELLIPSE 2 2 3 0 Black White
% 4 4254 1708 4312 1500 4254 1917 4254 1501
% 
{\color[named]{Black}{%
\special{pn 8}%
\special{pa 4254 1500}%
\special{pa 4258 1500}%
\special{pa 4258 1501}%
\special{pa 4261 1501}%
\special{fp}%
\special{pa 4280 1521}%
\special{pa 4280 1523}%
\special{pa 4281 1523}%
\special{pa 4281 1525}%
\special{pa 4282 1525}%
\special{pa 4282 1527}%
\special{fp}%
\special{pa 4294 1556}%
\special{pa 4294 1559}%
\special{pa 4295 1559}%
\special{pa 4295 1563}%
\special{fp}%
\special{pa 4303 1595}%
\special{pa 4303 1599}%
\special{pa 4304 1600}%
\special{pa 4304 1602}%
\special{fp}%
\special{pa 4309 1637}%
\special{pa 4309 1645}%
\special{fp}%
\special{pa 4312 1681}%
\special{pa 4312 1681}%
\special{pa 4312 1689}%
\special{fp}%
\special{pa 4312 1727}%
\special{pa 4312 1735}%
\special{fp}%
\special{pa 4309 1771}%
\special{pa 4309 1779}%
\special{fp}%
\special{pa 4304 1813}%
\special{pa 4304 1816}%
\special{pa 4303 1817}%
\special{pa 4303 1820}%
\special{fp}%
\special{pa 4296 1853}%
\special{pa 4295 1854}%
\special{pa 4295 1857}%
\special{pa 4294 1857}%
\special{pa 4294 1860}%
\special{fp}%
\special{pa 4283 1889}%
\special{pa 4283 1889}%
\special{pa 4282 1889}%
\special{pa 4282 1891}%
\special{pa 4281 1891}%
\special{pa 4281 1893}%
\special{pa 4280 1893}%
\special{pa 4280 1894}%
\special{fp}%
\special{pa 4261 1915}%
\special{pa 4258 1915}%
\special{pa 4258 1916}%
\special{pa 4254 1916}%
\special{fp}%
}}%
% ELLIPSE 2 0 3 0 Black White
% 4 4254 1708 4196 1500 4254 1501 4254 1917
% 
{\color[named]{Black}{%
\special{pn 8}%
\special{ar 4254 1708 58 208  1.5707963  4.7123890}%
}}%
% ELLIPSE 0 0 3 0 Black White
% 4 5896 4330 5802 3968 5896 3968 5896 4690
% 
{\color[named]{Black}{%
\special{pn 20}%
\special{ar 5896 4330 94 362  1.5707963  4.7123890}%
}}%
% ELLIPSE 0 2 3 0 Black White
% 4 5896 4330 5990 3968 5896 4690 5896 3968
% 
{\color[named]{Black}{%
\special{pn 20}%
\special{pa 5896 3968}%
\special{pa 5901 3968}%
\special{pa 5901 3969}%
\special{pa 5903 3969}%
\special{fp}%
\special{pa 5926 3986}%
\special{pa 5926 3987}%
\special{pa 5927 3988}%
\special{pa 5927 3989}%
\special{pa 5928 3989}%
\special{pa 5928 3990}%
\special{pa 5929 3991}%
\special{pa 5929 3992}%
\special{pa 5929 3992}%
\special{fp}%
\special{pa 5943 4018}%
\special{pa 5944 4018}%
\special{pa 5944 4020}%
\special{pa 5945 4020}%
\special{pa 5945 4022}%
\special{pa 5946 4022}%
\special{pa 5946 4023}%
\special{fp}%
\special{pa 5956 4051}%
\special{pa 5956 4053}%
\special{pa 5957 4053}%
\special{pa 5957 4056}%
\special{pa 5958 4056}%
\special{pa 5958 4057}%
\special{fp}%
\special{pa 5966 4089}%
\special{pa 5966 4091}%
\special{pa 5967 4091}%
\special{pa 5967 4095}%
\special{pa 5967 4095}%
\special{fp}%
\special{pa 5974 4129}%
\special{pa 5974 4131}%
\special{pa 5975 4131}%
\special{pa 5975 4136}%
\special{fp}%
\special{pa 5980 4170}%
\special{pa 5980 4171}%
\special{pa 5981 4172}%
\special{pa 5981 4178}%
\special{fp}%
\special{pa 5985 4214}%
\special{pa 5985 4219}%
\special{pa 5986 4219}%
\special{pa 5986 4221}%
\special{fp}%
\special{pa 5988 4257}%
\special{pa 5988 4264}%
\special{fp}%
\special{pa 5990 4301}%
\special{pa 5990 4309}%
\special{fp}%
\special{pa 5990 4347}%
\special{pa 5990 4355}%
\special{fp}%
\special{pa 5989 4392}%
\special{pa 5989 4394}%
\special{pa 5988 4395}%
\special{pa 5988 4399}%
\special{fp}%
\special{pa 5986 4435}%
\special{pa 5986 4440}%
\special{pa 5985 4441}%
\special{pa 5985 4443}%
\special{fp}%
\special{pa 5982 4478}%
\special{pa 5982 4480}%
\special{pa 5981 4481}%
\special{pa 5981 4485}%
\special{fp}%
\special{pa 5976 4521}%
\special{pa 5976 4523}%
\special{pa 5975 4523}%
\special{pa 5975 4527}%
\special{fp}%
\special{pa 5968 4561}%
\special{pa 5968 4565}%
\special{pa 5967 4565}%
\special{pa 5967 4568}%
\special{fp}%
\special{pa 5958 4601}%
\special{pa 5958 4601}%
\special{pa 5958 4603}%
\special{pa 5957 4604}%
\special{pa 5957 4607}%
\special{pa 5956 4607}%
\special{pa 5956 4607}%
\special{fp}%
\special{pa 5946 4638}%
\special{pa 5946 4638}%
\special{pa 5945 4638}%
\special{pa 5945 4640}%
\special{pa 5944 4640}%
\special{pa 5944 4642}%
\special{pa 5943 4642}%
\special{pa 5943 4643}%
\special{fp}%
\special{pa 5930 4668}%
\special{pa 5929 4668}%
\special{pa 5929 4670}%
\special{pa 5928 4670}%
\special{pa 5928 4671}%
\special{pa 5927 4671}%
\special{pa 5927 4672}%
\special{pa 5926 4673}%
\special{fp}%
\special{pa 5903 4691}%
\special{pa 5901 4691}%
\special{pa 5901 4692}%
\special{pa 5896 4692}%
\special{fp}%
}}%
% LINE 0 0 3 0 Black White
% 2 1694 4200 5894 3966
% 
{\color[named]{Black}{%
\special{pn 20}%
\special{pa 1694 4200}%
\special{pa 5894 3966}%
\special{fp}%
}}%
% DOT 0 0 3 0 Black White
% 2 5896 4692 5896 4692
% 
{\color[named]{Black}{%
\special{pn 4}%
\special{sh 1}%
\special{ar 5896 4692 16 16 0  6.28318530717959E+0000}%
\special{sh 1}%
\special{ar 5896 4692 16 16 0  6.28318530717959E+0000}%
}}%
% LINE 0 0 3 0 Black White
% 2 5896 4692 1696 4204
% 
{\color[named]{Black}{%
\special{pn 20}%
\special{pa 5896 4692}%
\special{pa 1696 4204}%
\special{fp}%
}}%
% DOT 2 0 3 0 Black White
% 2 1696 4204 1696 4204
% 
{\color[named]{Black}{%
\special{pn 4}%
\special{sh 1}%
\special{ar 1696 4204 8 8 0  6.28318530717959E+0000}%
\special{sh 1}%
\special{ar 1696 4204 8 8 0  6.28318530717959E+0000}%
}}%
% DOT 0 0 3 0 Black White
% 2 1696 4204 1696 4204
% 
{\color[named]{Black}{%
\special{pn 4}%
\special{sh 1}%
\special{ar 1696 4204 16 16 0  6.28318530717959E+0000}%
\special{sh 1}%
\special{ar 1696 4204 16 16 0  6.28318530717959E+0000}%
}}%
% ELLIPSE 2 0 3 0 Black White
% 4 3056 4244 3022 4124 3056 4124 3056 4362
% 
{\color[named]{Black}{%
\special{pn 8}%
\special{ar 3056 4244 34 120  1.5707963  4.7123890}%
}}%
% ELLIPSE 2 2 3 0 Black White
% 4 3056 4244 3090 4124 3056 4362 3056 4124
% 
{\color[named]{Black}{%
\special{pn 8}%
\special{pa 3056 4124}%
\special{pa 3059 4124}%
\special{pa 3059 4125}%
\special{pa 3061 4125}%
\special{pa 3061 4126}%
\special{pa 3062 4126}%
\special{fp}%
\special{pa 3078 4151}%
\special{pa 3078 4154}%
\special{pa 3079 4154}%
\special{pa 3079 4157}%
\special{pa 3080 4157}%
\special{pa 3080 4158}%
\special{fp}%
\special{pa 3087 4193}%
\special{pa 3087 4199}%
\special{pa 3088 4199}%
\special{pa 3088 4200}%
\special{fp}%
\special{pa 3090 4240}%
\special{pa 3090 4248}%
\special{fp}%
\special{pa 3088 4287}%
\special{pa 3088 4289}%
\special{pa 3087 4289}%
\special{pa 3087 4295}%
\special{fp}%
\special{pa 3080 4330}%
\special{pa 3080 4331}%
\special{pa 3079 4331}%
\special{pa 3079 4334}%
\special{pa 3078 4334}%
\special{pa 3078 4337}%
\special{fp}%
\special{pa 3062 4362}%
\special{pa 3061 4362}%
\special{pa 3061 4363}%
\special{pa 3059 4363}%
\special{pa 3059 4364}%
\special{pa 3056 4364}%
\special{fp}%
}}%
% ELLIPSE 2 2 3 0 Black White
% 4 4256 4280 4306 4058 4256 4418 4256 4058
% 
{\color[named]{Black}{%
\special{pn 8}%
\special{pa 4256 4058}%
\special{pa 4259 4058}%
\special{pa 4259 4059}%
\special{pa 4262 4059}%
\special{pa 4262 4060}%
\special{pa 4262 4060}%
\special{fp}%
\special{pa 4279 4084}%
\special{pa 4280 4084}%
\special{pa 4280 4086}%
\special{pa 4281 4087}%
\special{pa 4281 4089}%
\special{pa 4282 4089}%
\special{pa 4282 4090}%
\special{fp}%
\special{pa 4291 4121}%
\special{pa 4291 4123}%
\special{pa 4292 4124}%
\special{pa 4292 4128}%
\special{pa 4293 4129}%
\special{fp}%
\special{pa 4299 4163}%
\special{pa 4299 4163}%
\special{pa 4299 4170}%
\special{pa 4300 4171}%
\special{fp}%
\special{pa 4303 4208}%
\special{pa 4303 4210}%
\special{pa 4304 4211}%
\special{pa 4304 4215}%
\special{fp}%
\special{pa 4306 4253}%
\special{pa 4306 4261}%
\special{fp}%
\special{pa 4306 4300}%
\special{pa 4306 4308}%
\special{fp}%
\special{pa 4304 4345}%
\special{pa 4304 4349}%
\special{pa 4303 4350}%
\special{pa 4303 4353}%
\special{fp}%
\special{pa 4299 4389}%
\special{pa 4299 4389}%
\special{pa 4299 4397}%
\special{fp}%
\special{pa 4293 4430}%
\special{pa 4293 4432}%
\special{pa 4292 4432}%
\special{pa 4292 4436}%
\special{pa 4291 4436}%
\special{pa 4291 4436}%
\special{fp}%
\special{pa 4282 4469}%
\special{pa 4282 4471}%
\special{pa 4281 4471}%
\special{pa 4281 4473}%
\special{pa 4280 4474}%
\special{pa 4280 4475}%
\special{fp}%
\special{pa 4262 4500}%
\special{pa 4262 4500}%
\special{pa 4262 4501}%
\special{pa 4259 4501}%
\special{pa 4259 4502}%
\special{pa 4256 4502}%
\special{fp}%
}}%
% ELLIPSE 2 0 3 0 Black White
% 4 4256 4280 4206 4058 4256 4058 4256 4418
% 
{\color[named]{Black}{%
\special{pn 8}%
\special{ar 4256 4280 50 222  1.5707963  4.7123890}%
}}%
% LINE 2 0 3 0 Black White
% 2 5490 3200 2102 3200
% 
{\color[named]{Black}{%
\special{pn 8}%
\special{pa 5490 3200}%
\special{pa 2102 3200}%
\special{fp}%
}}%
% STR 2 0 3 0 Black White
% 4 1690 3100 1690 3200 5 0 0 0
% $t=0$
\put(16.9000,-32.0000){\makebox(0,0){$t=0$}}%
% STR 2 0 3 0 Black White
% 4 5890 3100 5890 3200 5 0 0 0
% $t=1$
\put(58.9000,-32.0000){\makebox(0,0){$t=1$}}%
% STR 2 0 3 0 Black White
% 4 3700 800 3700 900 5 0 0 0
% Figure of $\Delta $
\put(37.0000,-9.0000){\makebox(0,0){Figure of $\Delta $}}%
% STR 2 0 3 0 Black White
% 4 3700 3700 3700 3800 5 0 0 0
% Figure of $\Gamma  $
\put(37.0000,-38.0000){\makebox(0,0){Figure of $\Gamma  $}}%
\end{picture}%

\end{center}

%%%%%%%%%%%%%%%%%%%%%%%%%%%%%

\vspace{0.5cm}

\begin{lem}\label{L-lem0}
Suppose that $l$ is a prime number and $h^{20}(X)>0$.
Then $H_2(X,\Q)/\NF^B(X)$ has a basis 
$\{\sigma_*^i\Delta-\sigma_*^{i+1}\Delta\}_{0\leq i \leq l-2}$.
\end{lem}
\begin{pf}
Let $V\subset H_2(X,\Q)/\NF^B(X)$ be the subgroup generated by 
the image of $\{\sigma_*^i\Delta-\sigma_*^j\Delta\}_{i<j}$.
We want to show $V= H_2(X,\Q)/\NF^B(X)$.
Since $H_2(X,\Q)/\NF^B(X)$ is an irreducible $\Q[\sigma]$-module (by the proof of
Prop. \ref{NS=NF}) and $V$ is stable
under the action of $\Q[\sigma]$, it is 
enough to show $V\ne0$.
By Prop. \ref{goodbasis}, it is enough to show
\[
\int_{\Delta-\sigma_*\Delta}t^{i-1}dt\frac{dx}{y}=
(1-\zeta_l^i)\int_\Delta t^{i-1}dt\frac{dx}{y}\ne0.
\]
Due to {\bf (E1)},{\bf (E3)} and {\bf(E4)} there exist
3-distinct real roots $r_1(t),r_2(t),r_3(t)$
of $4x^3-g_2(t^l)x-g_3(t^l)$ for $0<t<1$.
Let them satisfy
$r_1(t)>r_2(t)>r_3(t)$ (resp. $r_1(t)<r_2(t)<r_3(t)$) if $\kappa>0$ (resp. $\kappa<0$). 
Then
\begin{equation}\label{Example-4}
\int_{\Delta} t^{j-1}dt\frac{dx}{y}=
2\sqrt{-1}
\overbrace{\int_0^1 t^{j-1}dt\int_{r_1(t)}^{r_2(t)}
\frac{dx}{\sqrt{-\kappa^{-1}(4x^3-g_2(t^l)x-g_3(t^l))}}}^{\R_{>0}}
\end{equation}
is not zero, so we are done.
\end{pf}
Similarly to \eqref{Example-4}, we have
\begin{equation}\label{Example-5}
\int_\Gamma t^{j-1}dt\frac{dx}{y}=
2\int_0^1 t^{j-1}dt\int_{r_2(t)}^{r_3(t)}
\frac{dx}{\sqrt{\kappa^{-1}(4x^3-g_2(t^l)x-g_3(t^l))}}\in\R_{>0}.
\end{equation}

\begin{lem}\label{L-lem1}
Let $F_\infty$ denotes the infinite Frobenius.
Then 
\[
F_\infty(\Delta)=-\Delta,\quad
F_\infty(\Gamma)=\Gamma.
\]
If $l$ is a prime number and $h^{20}(X)>0$,
then the fixed part $(H_2^B(X,\Q)/\NF^B(X))^{F_\infty=1}$ has a basis
\[
\sigma_*^i\Delta-\sigma_*^{-i}\Delta,\quad 1\leq i \leq \frac{l-1}{2}.
\]
\end{lem}
\begin{pf}
We show $F_\infty(\Delta)=-\Delta$.
Let $\delta_1\in H_1(f^{-1}(t),\Z)$ ($0<t<1$)
be the vanishing cycle as $t\to 1$.
Then it is enough to show $F_\infty(\delta_1)=-\delta_1$.
We keep the notation in the proof of Lemma \ref{L-lem0}.
Fix $0<t <1$ and $x\in [r_1(t),r_2(t)]$.  Then
$4x^3-g_1(t^l)x-g_3(t^l)\leq 0$ if $\kappa>0$ and $\geq 0$ if $\kappa<0$.
Therefore $y$ takes values in purely imaginary numbers, so that
$F_\infty(x,y)=(x,-y)$. This means $F_\infty\delta_1=-\delta_1$.
In the same way we have $F_\infty\delta_0=\delta_0$ where $\delta_0$
denotes the vanishing cycle as $t\to 0$. This implies 
$F_\infty(\Gamma)=\Gamma$ as well.
The last assertion follows from this and $F_\infty\sigma=\sigma^{-1}F_\infty$
together with Lemma \ref{L-lem0}.
\end{pf}
\subsection{Regulator indecomposable elements}

\begin{thm}\label{Example-thm0}
Suppose that $l$ is a prime number and $h^{20}(X)>0$.
Put $h=\dim F^1 H^2(X)_\ind$ and
$\zeta=\exp(2\pi i/l)$.
Let
\[
A=
\left((\zeta^{pq}-\zeta^{-pq})
\int_{\Delta} t^{p-1}dt\frac{dx}{y}\right)_{1\leq p\leq h,~1\leq q \leq (l-1)/2}
\]
be $h\times (l-1)/2$-matrix (the entries are real numbers by \eqref{Example-4}).
Then
\[
\mathrm{Ext}^1_{\R\mbox{-}\MHS}(\R,H^2(X)_\ind\ot\R(1))^{F_\infty=1}
\cong\Coker[A:\R^{(l-1)/2}\lra \R^h]
\]
and we have
\[
\reg_\R(\xi)=
\pm\left(
\int_{\Gamma} dt\frac{dx}{y},
%\int_{\Gamma} tdt\frac{dx}{y},
\cdots,
\int_{\Gamma} t^{h-1}dt\frac{dx}{y}
\right)\in \R^h/\Image A
\]
under the above isomorphism.
\end{thm}
\begin{pf}
The first assertion is obtained by applying Prop. \ref{goodbasis} and Lemma \ref{L-lem1}
to Theorem \ref{ExtMHS} (1).
The second assertion follows from Theorem \ref{ExtMHS} (2).
\end{pf}
\begin{cor}\label{Example-cor}
Suppose that $l$ is a prime number and $h^{20}(X)>0$. Then we have
\[
\reg_\R(\xi)\ne0\in\Ext^1_{\MHS}(\R,H^2(X)_\ind\ot\R(1))^{F_\infty=1}.
\]
In particular $\xi$ is real regulator indecomposable.
\end{cor}
\begin{pf}
Put
\[
I_p:=\int_\Delta t^{p-1}dt\frac{dx}{y},
\quad
J_p:=\int_\Gamma t^{p-1}dt\frac{dx}{y}.
\]
Then
\[
\reg_\R(\xi_D)\ne0\in\Ext^1_{\MHS}(\R,H^2(X)_\ind\ot\R(1))^{F_\infty=1}
\]
if and only if the rank of a matrix
\begin{equation}\label{Example-cor-1}
\begin{pmatrix}
 (\zeta-\zeta^{-1})I_{1}& (\zeta^2-\zeta^{-2})I_{1}
&\cdots&(\zeta^{(l-1)/2}-\zeta^{-(l-1)/2})I_{1}&J_1\\
 (\zeta^2-\zeta^{-2})I_{2}& (\zeta^{4}-\zeta^{-4})I_{2}
&\cdots&(\zeta^{l-1}-\zeta^{-(l-1)})I_{2}&J_2\\
 \vdots&\vdots&&\vdots&\vdots\\
(\zeta^{h}-\zeta^{-h})I_{h}& (\zeta^{2h}-\zeta^{-2h})I_{h}
&\cdots&(\zeta^{h(l-1)/2}-\zeta^{-h(l-1)/2})I_{h}&J_h\\
\end{pmatrix}
\end{equation}
is maximal.
It is enough to show that
\begin{equation}\label{Example-cor-2}
\det
\begin{pmatrix}
 (\zeta-\zeta^{-1})& (\zeta^2-\zeta^{-2})
&\cdots&(\zeta^{(l-1)/2}-\zeta^{-(l-1)/2})&J_1/I_1\\
 (\zeta^2-\zeta^{-2})& (\zeta^{4}-\zeta^{-4})
&\cdots&(\zeta^{l-1}-\zeta^{-(l-1)})&J_2/I_2\\
 \vdots&\vdots&&\vdots&\vdots\\
(\zeta^{k}-\zeta^{-k})& (\zeta^{2k}-\zeta^{-2k})
&\cdots&(\zeta^{k(l-1)/2}-\zeta^{-k(l-1)/2})&J_{k}/I_{k}\\
\end{pmatrix}
\end{equation}
is nonzero where $k=(l+1)/2$.
Since the sum of the $(k-1)$-th row and $k$-th row is
$
(0,\cdots,0,J_{k-1}/I_{k-1}+J_{k}/I_k)
$,
one has
\begin{align*}
\eqref{Example-cor-2}=&
(J_{k-1}/I_{k-1}+J_{k}/I_k)\times\det(\zeta^{pq}-\zeta^{-pq})_{1\leq p,q\leq (l-1)/2}\\
=&(J_{k-1}/I_{k-1}+J_{k}/I_k)\times\sqrt{(-l)^{(l-1)/2}}.
\end{align*}
Since $J_p/I_p\in i\R_{>0}$ by \eqref{Example-4} and \eqref{Example-5}, 
this is not zero. 
\end{pf}
%%%%%%%%%%%%%%%%%%%%%%%%%

\subsection{Explicit computation of regulator}
We show more on computation of real regulator in the following case
\[
X=X/\Q:-12y^2=4x^3-g_2(t^l)x-g_3(t^l),
\]
\[
(g_2,g_3)=(12(9-8t),-8(8t^{2}-36t+27))
\]
an elliptic surface defined over $\Q$ and
\[
\xi=\left[\frac{y-(x+1)}{y+(x+1)},D\right]\in H^3_{\cM}(X,\Q(2)).
\]
Here $D:=f^{-1}(1)$ is a multiplicative fiber of type $I_1$, and it splits over $\Q$.
One can show that if $l\geq 5$ is a prime number, then 
$\xi$ is integral in the sense of \cite{scholl},
namely, it lies in the image of the motivic cohomology group of a regular proper flat 
model ${\mathscr X}$ over $\Spec\Z$.
When $l=1$, $X$ is the universal elliptic curve over $X_1(3)$.
(However, if $l>1$ it is no longer a universal elliptic curve for congruence subgroup.)
Let $q=\exp(2\pi i z)$ and
\[
E_{3a}(z):=1-9\sum_{n=1}^\infty
\left(\sum_{k|n}\left(\frac{k}{3}\right)k^2\right)q^n,
\]
\[
E_{3b}(z):=\sum_{n=1}^\infty\left(\sum_{k|n}\left(\frac{n/k}{3}\right)k^2\right)q^n
\]
be the Eisenstein series of weight 3 for $\Gamma_1(3)$, where $(\frac{k}{3})$ denotes 
the Legendre symbol.
Then
\[
t^l
=\frac{E_{3a}}{E_{3a}+27E_{3b}}
\]
and
\[
l\frac{dt}{t}\frac{dx}{y}=-27E_{3b}\frac{du}{u}\frac{dq}{q},\quad
\frac{lt^{l-1}dt}{t^l-1}\frac{dx}{y}=E_{3a}\frac{du}{u}\frac{dq}{q}
\]
where ``$du/u$" denotes the canonical invariant 1-form of the Tate curve around 
the cusp $z=i\infty$ ($t=1$).
Therefore we have
\begin{equation}\label{Example-1}
\int_{\Delta} t^{j-1}dt\frac{dx}{y}=\frac{-27}{l}\times (2\pi i)^2
\int_0^{i\infty} t^{j}E_{3b}(z)dz
\end{equation}
\begin{equation}\label{Example-2}
\int_\Gamma t^{j-1}dt\frac{dx}{y}=\frac{-27}{l}\times (2\pi i)^2
\int_0^{i\infty} t^{j}E_{3b}(z)zdz.
\end{equation}
On the other hand there are formulas
\begin{equation}\label{Example-3}
\frac{E_{3a}}{E_{3a}+27E_{3b}}(\frac{-1}{3z})
=\frac{27E_{3b}}{E_{3a}+27E_{3b}}(z),\quad
27E_{3b}(\frac{-1}{3z})=
3\sqrt{3}iz^3 E_{3a}(z)
\end{equation}
on the Eisenstein series.
Applying \eqref{Example-3} to \eqref{Example-1} and \eqref{Example-2},
we have the following theorem.
%%%%%%%%%%%%%%%%%%%%%%%%%%%%

\begin{thm}\label{Example-thm}
Put $c:=\exp(-2\pi/\sqrt{3})=0.026579933\cdots$. 
Define rational numbers $a_n(j)$ and $b_n(j)$ by
\begin{multline*}
E_{3b}\left(\frac{E_{3a}}{E_{3a}+27E_{3b}}\right)^{j/l}=\sum_{n=1}^\infty a_n{(j)}q^n\\
=q+\left(3-27\frac{j}{l}\right)q^2
+\left(9-\frac{81}{2}\frac{j}{l}+\frac{729}{2}\left(\frac{j}{l}\right)^2\right)q^3+\cdots,
\end{multline*}
\begin{multline*}
E_{3a}\left(\frac{E_{3b}}{q(E_{3a}+27E_{3b})}\right)^{j/l}
=\sum_{n=0}^\infty b_n{(j)}q^n\\
=1+\left(-9-15\frac{j}{l}\right)q+
\left(27+\frac{387}{2}\frac{j}{l}+\frac{225}{2}\left(\frac{j}{l}\right)^2\right)q^2+\cdots.
\end{multline*}
Put
\begin{align*}
I(j)&=\sum_{n=1}^\infty \frac{a_n{(j)}}{n}c^n+3^{3j/l-3}\sum_{n=0}^\infty b_n{(j)}\left(\frac{1}{n+j/l}
+\frac{\sqrt{3}}{2\pi(n+j/l)^2}\right)c^{n+j/l}\\
J(j)&=\sum_{n=1}^\infty a_n(j)\left(
\frac{2\pi}{\sqrt{3}n}+\frac{1}{n^2}\right)c^n+2\pi\cdot 3^{3j/l-7/2}\sum_{n=0}^\infty 
\frac{b_n(j)}{n+j/l}c^{n+j/l}.
\end{align*}
Then we have
\[
\int_{\Delta} t^{j-1}dt\frac{dx}{y}
=\frac{54\pi i}{l}I(j),\quad 
\int_\Gamma t^{j-1}dt\frac{dx}{y}
=\frac{-27}{l}J(j)
\]
for $1\leq j\leq l-1$.
\end{thm}
This is useful since the series $I(j)$ and $J(j)$ converge rapidly !
%}\end{exmp}
%%%%%%%%%%%%%%%%%%%%%%%%%%%%%

\begin{exmp}\label{exp-l5}{\rm Suppose $l=5$. Then $X$ is a K3 surface.
By Thm.\ref{Example-thm}, one has 
\begin{center}
\begin{tabular}{|c|c|c|}
\hline
&$I(j)$
&$J(j)$\\
\hline
$j=1$&
$ 0.42745977255318$&$0.717696894965804$
\\
\hline
$j=2$&
$0.151180954233147$
&$0.377159120670032$
\\
\hline
$j=3$&
$0.0871841692346256$
&$0.261572572611421$
\\
\hline
$j=4$&
$0.0603840144077692$&$0.202670503662525$
\\
\hline
\end{tabular}
\end{center}
\[
\mathrm{Ext}^1_{\R\mbox{-}\MHS}(\R,H^2(X)_\ind\ot\R(1))^{F_\infty=1}
\cong\Coker(\R^2\os{A}{\lra} \R^3).
\]
Since this is 1-dimensional,
this has the canonical base $e_{\ind,\Q}$ (up to $\Q^\times$) and a {\it different} base
$e_{\ind,\Q}^\ff$
(\S \ref{false-sect}).
With respect to $e_{\ind,\Q}^\ff$, one has
\[
\reg_\R(\xi_D)=\pi^2
\begin{vmatrix}
 i(\zeta-\zeta^{-1})I(1)& i(\zeta^2-\zeta^{-2})I(1)&J(1)\\
 i(\zeta^2-\zeta^{-2})I(2)& i(\zeta^4-\zeta^{-4})I(2)&J(2)\\
i(\zeta^3-\zeta^{-3})I(3)& i(\zeta^6-\zeta^{-6})I(3)&J(3)\\
\end{vmatrix}\mod \Q^\times \quad(\zeta:=\exp(2\pi i/5)).
\]
Since $s=(l-1)/2=2$ and $\det H_\dR^2(X/\Q)_\ind\ot
[\det H_B^2(X_\C)_\ind]^{-1}=\sqrt{5}$, %(see Prop.\ref{qstr-8} for the notations),
one has $e_{\ind,\Q}^\ff=(2\pi\sqrt{-1})^{-2}\sqrt{5}~e_{\ind,\Q}
=\sqrt{5}\pi^{-2}~e_{\ind,\Q}$ mod $\Q^\times$
by Prop.\ref{qstr-8}.  Hence
\begin{align*}
\reg_\R(\xi_D)
&=\frac{\sqrt{5}}{\pi^2}\cdot\pi^2
\begin{vmatrix}
 i(\zeta-\zeta^{-1})I(1)& i(\zeta^2-\zeta^{-2})I(1)&J(1)\\
 i(\zeta^2-\zeta^{-2})I(2)& i(\zeta^4-\zeta^{-4})I(2)&J(2)\\
i(\zeta^3-\zeta^{-3})I(3)& i(\zeta^6-\zeta^{-6})I(3)&J(3)\\
\end{vmatrix}\\
&=-5\sqrt{5}I(1)I(2)I(3)\left(\frac{J(2)}{I(2)}+\frac{J(3)}{I(3)}\right)\\
&=0.346139631939354
\mod \Q^\times 
\end{align*}
with respect to $e_{\ind,\Q}$.
}
\end{exmp}

\begin{exmp}\label{exp-l7} {\rm Suppose $l=7$. Then $h^{20}(X)=h^{02}(X)=2$, $h^{11}(X)=30$.
\begin{center}
\begin{tabular}{|c|c|c|}
\hline
&$I(j)$
&$J(j)$\\
\hline
$j=1$&
$0.740059830730164$&$0.987994510350351$
\\
\hline
$j=2$&
$0.24646699651114$
&$0.51401702238944$
\\
\hline
$j=3$&
$0.137265313181901$
&$0.354195498081428$
\\
\hline
$j=4$&
$0.0929578147374374$&$0.273237679671921$
\\
\hline
$j=5$&
$0.0696363855176379$&$0.224004116344261$
\\
\hline
$j=6$&
$0.0554349861351089$&$0.19073921727221$
\\
\hline
\end{tabular}
\end{center}
\[
\mathrm{Ext}^1_{\R\mbox{-}\mathrm{MHS}}(\R,H^2(X)_\ind\ot\R(1))^{F_\infty=1}
\cong\Coker(\R^3\os{A}{\lra} \R^4).
\]
Since $s=(l-1)/2=3$ and $\det H_\dR^2(X/\Q)_\ind\ot
[\det H_B^2(X_\C)_\ind]^{-1}=\sqrt{-7}$,
one has $e_{\ind,\Q}^\ff=(2\pi\sqrt{-1})^{-3}\sqrt{-7}~e_{\ind,\Q}
=\sqrt{7}\pi^{-3}~e_{\ind,\Q}$ mod $\Q^\times$, and
\begin{align*}
\reg_\R(\xi_D)&=\frac{\sqrt{7}}{\pi^3}\cdot\pi^3
\begin{vmatrix}
 i(\zeta-\zeta^{-1})I(1)& i(\zeta^2-\zeta^{-2})I(1)& i(\zeta^3-\zeta^{-3})I(1)&J(1)\\
 i(\zeta^2-\zeta^{-2})I(2)& i(\zeta^4-\zeta^{-4})I(2)& i(\zeta^6-\zeta^{-6})I(2)&J(2)\\
i(\zeta^3-\zeta^{-3})I(3)& i(\zeta^6-\zeta^{-6})I(3)& i(\zeta^9-\zeta^{-9})I(3)&J(3)\\
i(\zeta^4-\zeta^{-4})I(4)& i(\zeta^8-\zeta^{-8})I(4)& i(\zeta^{12}-\zeta^{-12})I(4)&J(4)\\
\end{vmatrix}\\
&=49I(1)I(2)I(3)I(4)\left(\frac{J(3)}{I(3)}+\frac{J(4)}{I(4)}\right)\\
&=0.629487860860585
\mod \Q^\times\quad(\zeta:=\exp(2\pi i/7))
\end{align*}
with respect to the canonical $\Q$-structure $e_{\ind,\Q}$.
}
\end{exmp}
\begin{rem}
According to the Beilinson conjecture, $\reg_\R(\xi_D)$
in Example \ref{exp-l5} or \ref{exp-l7} is expected to be the value of 
the $L$-function $L'(h^2(X)_\ind,1)$ $(\cite{schneider})$. 
\end{rem}
%%%%%%%%%%%%%%%%%%%%%%%%%%%%%%%%%%%%%%%%%%%%%%%

\def\EE{2g_2g'_3-3g'_2g_3}
\def\EEE{6g_2g'_3-9g'_2g_3}

\section{Appendix : Gauss-Manin connection for a hyperelliptic fibration}\label{Appendix-sect}
We work over a field $K$ of characteristic zero.
For a smooth scheme $Y$ over $T$, we denote by $\Omega^q_{Y/T}=\os{q}{\wedge}_{\O_Y}
\Omega^1_{Y/T}$ the sheaf of relative differential $q$-forms on $Y$ over $T$. 
If $T=\Spec K$, we
simply write $\Omega^q_{Y}=\Omega^q_{Y/K}$.

\medskip

In this section, we discuss the {\it Gauss-Manin connection}
\[
\nabla:R^1f_*\Omega^\bullet_{U/S}\lra \Omega^1_S\ot
R^1f_*\Omega^\bullet_{U/S}%\quad (\Omega^\bullet_S:=\Omega^\bullet_{S/R})
\]
for $U/S$ a smooth proper family of hyperelliptic curves.
This is defined to be the connecting homomorphism
\begin{equation}\label{gm00}
R^1f_*\Omega^\bullet_{U/S}\to R^2f_*(f^*\Omega^1_S\ot
\Omega^{\bullet-1}_{U/S})\cong \Omega^1_S\ot R^2f_*(
\Omega^{\bullet-1}_{U/S})\cong \Omega^1_S\ot R^1f_*\Omega^\bullet_{U/S}
\end{equation}
which arises from an exact sequence
\[
0\lra
f^*\Omega^1_S\ot\Omega^{\bullet-1}_{U/S}
\lra
\bar{\Omega}_U^\bullet
\lra \Omega^\bullet_{U/S}\lra 0,\quad 
\bar{\Omega}^\bullet_U:=\Omega^\bullet_U/\Image(f^*\Omega^2_S\ot\Omega^{\bullet-2}_U)
\]
(cf. \cite{h} Ch.III, \S 4).
Here the first isomorphism in \eqref{gm00} is the projection formula,
and the second one is due to the identification 
$R^qf_*\Omega^{\bullet-1}_{U/S}\cong R^{q-1}f_*\Omega^{\bullet}_{U/S}$ with
which we should be careful about ``sign".
Indeed
the differential of the complex $\Omega^{\bullet-1}_{U/S}$ is ``$-d$"
\[
\Omega^{\bullet-1}_{U/S}:\O_U\os{-d}{\lra} \Omega^1_{U/S}
\quad\mbox{(the first term is placed in degree 1)}
\]
so that we need to arrange
 the sign to make an isomorphism between $R^qf_*\Omega^{\bullet}_{U/S}$
and $R^{q+1}f_*\Omega^{\bullet-1}_{U/S}$.
We make it by a commutative diagram 
\begin{equation}\label{gm-diag}
\xymatrix{
\O_U\ar[r]^d\ar[d]_{-\id}&\Omega^1_{U/S}\ar[d]^{\id}\\
\O_U\ar[r]^{-d}&\Omega^1_{U/S}
}
\end{equation}
Then $\nabla$ satisfies the usual Leibniz rule
\[
\nabla(gx)=dg\ot x+g\nabla(x),\quad x\in \vg(S,R^1f_*\Omega^\bullet_{U/S}),~g\in\O(S).
\]

\subsection{Family of hyperelliptic curves}\label{expPF-sect}
Let $S$ be an irreducible affine smooth variety over $K$.
Let $f(x)\in \O_S(S)[x]$ be a polynomial of degree $2g+1$ or $2g+2$
which has no multiple roots over any geometric points
$\bar{x}\in S$.
Then it defines a smooth family of hyperelliptic curves $f:U\to S$ 
defined by the Weierstrass equation $y^2=f(x)$. 
To be more precise, let $z=1/x,~u=y/x^{g+1}$ and put $g(z)=z^{2g+2}f(1/z)$.
Let
\[
U_0=\Spec \O_S(S)[x,y]/(y^2-f(x)),\quad
U_\infty=\Spec \O_S(S)[z,u]/(u^2-g(z)).
\]
Then $U$ is obtained by gluing
$U_0$ and $U_\infty$ via identification 
$z=1/x,~u=y/x^{g+1}$.
We assume that there is a section $e:S\to U$.
\begin{equation}\label{he0}
x^{i-1}\frac{dx}{y}=-z^{g-i}\frac{dz}{u},\quad 
\frac{y}{x^i}=\frac{u}{z^{g+1-i}},\quad 1\leq i \leq g.
\end{equation}
We shall compute the Gauss-Manin connection
\begin{equation}\label{gm4}
\nabla:H^1_\dR(U/S)\lra \Omega^1_S\otimes H^1_\dR(U/S),
\quad H^q_\dR(U/S):=H^q_\zar(U,\Omega^\bullet_{U/S})
\end{equation}
(we use the same symbol ``$\Omega^1_S$" for $\vg(S,\Omega^1_S)$ 
since it will be clear from the context which is meant).
To do this, we describe the de Rham cohomology
in terms of the Cech complex.
Write
\[
\check{C}^0({\mathscr F}):=\vg(U_0,{\mathscr F})\op\vg(U_\infty,{\mathscr F}),\quad
\check{C}^1({\mathscr F}):=\vg(U_0\cap U_\infty,{\mathscr F})
\]
for a (Zariski) sheaf $\mathscr F$.
Then the double complex
\[
\xymatrix{
\check{C}^0(\O_U)\ar[r]^d\ar[d]_\delta&
\check{C}^0(\Omega^1_{U/S})\ar[d]^\delta&(x_0,x_\infty)\ar[d]^\delta
\ar[r]^d&(dx_0,dx_\infty)\\
\check{C}^1(\O_U)\ar[r]^d&
\check{C}^1(\Omega^1_{U/S})&x_0-x_\infty
}
\]
gives rise to the total complex
\[
\check{C}^\bullet(U/S):\check{C}^0(\O_U)\os{\delta\times d}{\lra} 
\check{C}^1(\O_U)\times \check{C}^0(\Omega^1_{U/S})
\os{(-d)\times\delta}{\lra} \check{C}^1(\Omega^1_{U/S})
\]
of $R$-modules starting from degree 0, and the cohomology of it is the de Rham cohomology
$H^\bullet_\dR(U/S)$:
\[
H^q_\dR(U/S)=H^q(\check{C}^\bullet(U/S)),\quad q\geq 0.
\]
Elements of $H^1_\dR(U/S)$ are represented by cocycles
\[
(f)\times (x_0,x_\infty)\quad \mbox{with }df=x_0-x_\infty.
\]
\begin{lem}\label{he-basis}
Suppose
\[
f(x)=a_0+a_1x+\cdots+a_nx^n,\quad a_i\in \O(S)
\]
with $n=2g+1$ or $2g+2$. Put
\begin{equation}\label{he1}
\omega_i:=(0)\times\left(\frac{x^{i-1}dx}{y},-\frac{z^{g-i}dz}{u}\right),
\end{equation}
\begin{equation}\label{he2}
\omega_i^*:=\left(\frac{y}{x^i}\right)\times
\left(\left(\sum_{m>i}(m/2-i)a_mx^{m-i-1}\right)\frac{dx}{y},
\left(\sum_{m\leq i}(m/2-i)a_mz^{g-m+i}\right)\frac{dz}{u}\right)
\end{equation}
for $1\leq i \leq g$.
Then they give a basis of $H^1_\dR(U/S)$.
Moreover \eqref{he1} span the image of $\vg(U,\Omega^1_{U/S})\hra H^1_\dR(U/S)$.
\end{lem}
\begin{pf}
Exercise.
\end{pf}
\begin{lem}\label{cech-equiv}
There are the following equivalence relations.
\[
(x^iy^j)\times(0,0)\equiv(0)\times(-d(x^iy^j),0)\mod \Image\check{C}^0(\O_U),
\]
\[
(z^iu^j)\times(0,0)\equiv(0)\times(0,d(z^iu^j))\mod \Image\check{C}^0(\O_U).
\]
\end{lem}
\begin{pf}
Straightforward from the definition .
\end{pf}

\subsection{Computation of Gauss-Manin connection}
Let us compute $\nabla(\omega_i)$ and $\nabla(\omega_i^*)$.
Recall that there is the exact sequence
\begin{equation}\label{gm0}
0\lra \check{C}^\bullet(f^*\Omega^1_S\ot\Omega^{\bullet-1}_{U/S})\lra
\check{C}^\bullet(\bar{\Omega}_U^\bullet)
\lra \check{C}^\bullet(\Omega^\bullet_{U/S})\lra 0
\end{equation}
and it gives rise to the connecting homomorphism 
\[
%\begin{equation}\label{gm1}
\delta:H^1_\dR(U/S)=H^1(\check{C}^\bullet(\Omega^\bullet_{U/S}))\lra H^2(U,
f^*\Omega^1_S\ot\Omega^{\bullet-1}_{U/S})
=H^2(\check{C}^\bullet(f^*\Omega^1_S\ot\Omega^\bullet_{U/S})).
\]
Recall the isomorphism
\[
\check{C}^\bullet(f^*\Omega^1_S\ot\Omega^{\bullet-1}_{U/S})
\os{\cong}{\lra}
\Omega^1_S\ot\check{C}^\bullet(\Omega^{\bullet}_{U/S})
\]
induced from \eqref{gm-diag}. It induces the isomorphism 
\[
%\begin{equation}\label{gm2}
\iota:H^2(U,f^*\Omega^1_S\ot\Omega^{\bullet-1}_{U/S})\os{\cong}{\lra}
\Omega^1_S\ot H^1_\dR(U/S).
\]
By definition we have $\nabla=\iota\delta$ 
the Gauss-Manin connection \eqref{gm4}.
Let us write down the maps $\delta$ and $\iota$ in terms of Cech cocycles.
The differential operator $\cD$ on the total complex of the middle term of \eqref{gm0} is 
given as follows
\[\cD:
\check{C}^1(\O_U)\times \check{C}^0(\bar{\Omega}^1_{U})
\lra \check{C}^1(\bar{\Omega}^1_{U})\times \check{C}^0(\bar{\Omega}^2_{U}),
\]
\[
(\alpha)\times(\beta_0,\beta_\infty)\longmapsto
(-d\alpha+\beta_0-\beta_\infty)\times(d\beta_0,d\beta_\infty).
\]
We denote a lifting of $(z_0,z_\infty)\in \check{C}^0({\Omega}^1_{U/S})$
by $(\hat{z}_0,\hat{z}_\infty)\in \check{C}^0(\bar{\Omega}^1_U)$.
Then for $(\alpha)\times (z_0,z_\infty)\in H^1_\dR(U/S)$ one has
\begin{align}
(\alpha)\times (z_0,z_\infty)&\os{\delta}{\longmapsto} 
\cD((\alpha)\times (\hat{z}_0,\hat{z}_\infty))\\
&=(-d\alpha+\hat{z}_0-\hat{z}_\infty)\times(d\hat{z}_0,d\hat{z}_\infty)\\
&\in \check{C}^1(f^*\Omega^1_S)\times \check{C}^0(f^*\Omega^1_S\ot\Omega^1_{U/S}).
\end{align}
The isomorphism $\iota$ is given by
\begin{equation}\label{gm3}
(gdt)\times(dt\wedge z_0,dt \wedge z_\infty )
\os{\iota}{\longmapsto}
dt\ot [(-g)\times(z_0,z_\infty)]
\end{equation}
(the ``sign" appears in the above due to \eqref{gm-diag}).

To compute $\nabla(\omega_i)$ and $\nabla(\omega_i^*)$ for the basis
in Lemma \ref{he-basis}, there remains to compute lifting of $dx/y$ and $dz/u$.

\begin{lem}\label{ABCD}
Let $A,B\in\O(S)[x]$ and $C,D\in \O(S)[z]$ satisfy
\[
Af+B\frac{\partial f}{\partial x}=1,\quad Cg+D\frac{\partial g}{\partial z}=1.
\]
Put differential 1-forms
\[
\widehat{\frac{dx}{y}}:=\frac{Afdx+Bdf}{y}=Aydx+2Bdy
\in\vg(U_0,\Omega^1_U),
\]
\[
\widehat{\frac{dz}{u}}:=\frac{Cgdz+Ddg}{u}=Cudz+2Ddu\in\vg(U_\infty,\Omega^1_U).
\]
Then
\[
x^i\widehat{\frac{dx}{y}}\in\vg(U_0,\Omega^1_U),\quad
z^i\widehat{\frac{dz}{u}}\in\vg(U_\infty,\Omega^1_U)
\]
are liftings of 
$x^idx/y\in\vg(U_0,\Omega^1_{U/S})$ and
$z^idz/u\in\vg(U_\infty,\Omega^1_{U/S})$ respectively.
\end{lem}
\begin{pf}
Straightforward.
\end{pf}

By using the liftings in Lemma \ref{ABCD}, one can compute 
the map $\delta$.
With use of Lemma \ref{cech-equiv}, 
one finally obtains the connection matrix of $\nabla$.

\medskip

Here is an explicit formula in case of elliptic fibration (the proof is left to the reader).

\begin{thm}\label{A1}
Let $S$ be a smooth affine curve and $f:U\to S$ a projective smooth family of
elliptic curves whose affine form is given by a Weierstrass equation
$y^2=4x^3-g_2x-g_3$ with $\Delta:=g_2^3-27g_3^2\in \O_S(S)^\times$.
Suppose that $\Omega^1_S$ is a free $\O_S$-module with a base 
$dt\in \vg(S,\Omega^1_{S})$.
For $f\in \O_S(S)$, we define $f'\in\O_S(S)$ by $df=f'dt$.
Let
\begin{equation}\label{A1-1}
\canh:=(0)\times (\frac{dx}{y},-\frac{dz}{u})%\in \vg(U,\Omega^1_{U/S})
\end{equation}
\begin{equation}\label{A1-2}
\can:=(\frac{y}{x})\times (\frac{2xdx}{y},\frac{(g_2z+2g_3z^2)dz}{2u})
\end{equation}
be elements in $H^1_\dR(U/S)$. 
Then we have
\begin{equation}\label{A1-3}
\nabla\left(\canh\right)=\left(\frac{3(\EE)}{4\Delta}dt\ot\can
-\frac{\Delta^\prime}{12\Delta}dt\ot\canh\right)\in \Omega^1_S \ot H^1_\dR(U/S),
\end{equation}
\begin{equation}\label{A1-4}
\nabla\left(\can\right)=\left(\frac{\Delta^\prime}{12\Delta}dt\ot\can
-\frac{g_2(\EE)}{4\Delta}dt\ot\canh\right)\in 
\Omega^1_S \ot H^1_\dR(U/S).
\end{equation}
\end{thm}

\subsection{Deligne's canonical extension and the limiting Hodge filtration}\label{de-sect} 
Let $S$ be a smooth curve over $\C$
and $(\cH,\nabla)$ a vector bundle with integrable connection
over $S^*:=S-\{P\}$. Let $j:S^*\hra S$.
Then there is unique subbundle $\cH_e\subset j_*\cH$ satisfying the following conditions
(cf. \cite{zucker} (17)).
\begin{itemize}
\item
The connection extends to have log pole, 
$\nabla:\cH_e\to\Omega^1_S(\log P)\ot\cH_e$,
\item
each eigenvalue $\alpha$ of $\Res_P(\nabla)$ satisfies $0\leq \mathrm{Re}(\alpha)<1$.
\end{itemize}
The extended bundle $(\cH_e,\nabla)$ is called {\it Deligne's canonical extension}.
The inclusion map 
\[
[\cH_e\os{\nabla}{\to}\Omega^1_S(\log P)\ot\cH_e]\lra
[j_*\cH\os{\nabla}{\to}\Omega^1_{S^*}\ot j_*\cH]
\]
is a quasi-isomorphism of complexes of sheaves. 
%\begin{pf}
%Let $t\in \O_{S,P}$ be a uniformizer at $P$, and put ${\mathscr F}_{i}:=t^{-i}\cH_e$.
%Then $\nabla({\mathscr F}_{i})\subset \Omega^1_S(\log P)\ot {\mathscr F}_{i}$.
%Therefore it is enough to show that $\Res^{(i)}_P:\Gr^{{\mathscr F}}_i\cH_e\to 
%\Omega^1_S(\log P)\ot \Gr^{{\mathscr F}}_i\cH_e\cong\Gr^{{\mathscr F}}_i\cH_e$ is bijective
%unless $i=0$. This follows from the fact that the eigenvalues of $\Res^{(i)}_P$
%are in $[-i,-i+1)$.
%\end{pf}
Besides $\exp(-2\pi i\Res_P(\nabla))$ coincides with the monodromy
operator on $H_\C=\ker(\nabla^{\mathit an})$ around $P$
(cf. \cite{steenbrink}, (2.21)).

Let $(H_\Z,\cH,F^\bullet,\nabla)$ be a polarized VHS on $S$.
Then the eigenvalues of $\Res_P(\nabla)$ are in $\Q$ (\cite{Schmid} (4.5)). 
Moreover by the nilpotent orbit theorem (\cite{Schmid} (4.9)), 
one can show that $\hat{F}^\bullet:=\cH_e\cap j_*F^\bullet$ are subbundles of $\cH_e$
(e.g. \cite{msaito1} (2.2)).
$\hat{F}^\bullet$ and the $V$-filtration define the 
{\it limiting Hodge filtration} on $\cH_e\ot\C_p$ 
(note that $\hat{F}^\bullet\ot\C_P$ does {\it not} necessarily coincide with
the limiting
Hodge filtration unless the monodromy is unipotent.
See \cite{msaito2} (3.5) for the detail).

If $(H_\Z,\cH,F^\bullet,\nabla)$ is a VHS arising from a projective flat
family $f:X\to S$ such that $f$ is smooth over $S^*$
and $D_\mathrm{red}:=(f^{-1}(0))_\mathrm{red}$ is a NCD, then one has
\begin{equation}\label{g-de0}
\cH_e\cong R^qf_*\Omega^\bullet_{X/S}(\log D), \quad 
\hat{F}^i\cong R^qf_*\Omega^{\bullet\geq i}_{X/S}(\log D)
\end{equation}
(\cite{zucker} p.130, Corollary).
Put $\ol{\cH}_e:=\Coker[\Res_P(\nabla):\cH_e\ot\C_P\to \cH_e\ot\C_P]$ and
let 
\begin{equation}\label{residue}
\Res_p^{\cH_e}:H^1_\dR(S,\cH_e)=H^1_\zar(S,\cH_e\to\Omega^1_S(\log P)\ot\cH_e)
\lra \ol{\cH}_e
\end{equation}
be the map induced from a commutative diagram
\[
\xymatrix{
\cH_e\ar[r]\ar[d]&0\ar[d]\\
\Omega^1_S(\log P)\ot\cH_e\ar[r]^{\qquad\Res_P}&\ol{\cH}_e.
}
\]
If $(H_\Z,\cH,F^\bullet,\nabla)$ is the case \eqref{g-de0}, then
\eqref{residue} is compatible with the residue map
\[
\Res_D:H^{q+1}_\dR(X-D)\lra H^\dR_{2\dim X-q-2}(D)
\]
under the natural maps $H^1_\dR(S,\cH_e)\to H^{q+1}_\dR(X-D)$ and
$\ol{\cH}_e\to H^\dR_{2\dim X-q-2}(D)$.

\medskip

We have seen how to compute a connection matrix of the Gauss-Manin connection
for a family of hyperelliptic curves.
Once we have it, we can get $(\cH_e,\hat{F}^\bullet)$ automatically.
In case of an elliptic fibration, they are simply given as follows.
\begin{thm}\label{A2}
Let $f:U\to S:=\Spec\C[[t]]$ be an elliptic fibration defined by a minimal
Weierstrass equation $y^2=4x^3-g_2x-g_3$ with $g_2,g_3\in \C[[t]]$, $\Delta:=g_2^3-27g_3^2
\ne0$.
Then $\cH_e$ has a basis $\{\omega,\omega^*\}$ (resp. $\{t\omega,\omega^*\}$)
if $f$ has a semistable or smooth reduction (resp. additive reduction).
\end{thm}
Since we have Theorem \ref{A1},
we can show the above by case-by-case analysis based on $(\ord(g_2),\ord(g_3))$.
The detail is left to the reader.

%%%%%%%%%%%

\bigskip

\noindent
Department of Mathematics, Hokkaido University,
Sapporo 060-0810,
JAPAN

\medskip

\noindent
asakura@math.sci.hokudai.ac.jp

\end{document}